\documentclass[10pt]{amsart}

\usepackage{multirow}
\usepackage{graphicx}
\usepackage[top=26mm, bottom=26mm, left=23mm, right=25mm]{geometry}
\usepackage[dvips]{color}
\usepackage{bm}
\theoremstyle{definition}
\newtheorem{remark}{Remark}
\numberwithin{remark}{section}

\numberwithin{definition}{section}

\input{algorithm_labels.aux}

\newcommand{\vct}[1]{\bm{#1}}
\newcommand{\mtx}[1]{\mathsf{#1}}

\theoremstyle{definition}

\newcommand{\lsp}{\vspace{3mm}}

\newcommand{\vtwo}[2]{\left[\begin{array}{c} #1 \\ #2 \end{array}\right]}

\newcommand{\mtwo}[4]{\left[\begin{array}{cc}    #1 & #2 \\ #3 & #4  \end{array}\right]}

\begin{document}

\begin{center}
\textbf{\large An $O(N)$ algorithm for constructing the solution operator to
2D elliptic boundary value problems in the absence of body loads}

\lsp

{\small A. Gillman$^1$, and P.G. Martinsson$^2$\\
$1$ Department of Mathematics, Dartmouth College, $2$ Department of Applied Mathematics, University of Colorado at Boulder}

\lsp

\begin{minipage}{135mm}
\noindent\textbf{Abstract:}  The large sparse linear systems arising from the finite element or finite difference
discretization of elliptic PDEs can be solved directly via, e.g., nested dissection
or multifrontal methods.
Such techniques reorder the nodes in the grid to reduce
the asymptotic complexity of Gaussian elimination from $O(N^{2})$ to $O(N^{1.5})$
for typical problems in two dimensions.
It has recently been demonstrated that the complexity
can be further reduced to $O(N)$ by exploiting structure in the dense matrices that
arise in such computations (using, e.g., $\mathcal{H}$-matrix arithmetic).
This paper
demonstrates that such \textit{accelerated} nested dissection techniques become particularly effective for boundary value
problems without body loads when the solution is sought for several different sets
of boundary data, and the solution is required only near the boundary (as happens, e.g.,
in the computational modeling of scattering problems, or in engineering design of
linearly elastic solids).
In this case, a modified version of the accelerated nested dissection scheme
can execute any solve beyond the first in $O(N_{\rm boundary})$ operations,
where $N_{\rm  boundary}$ denotes the number of points on the boundary.
Typically, $N_{\rm boundary} \sim N^{0.5}$.
Numerical examples demonstrate the effectiveness of the procedure for a broad range of
elliptic PDEs that includes both the Laplace and Helmholtz equations.
\end{minipage}
\end{center}

\section{Introduction}
\label{sec:intro}

\subsection{Problem formulation}
This paper presents a fast solver for homogeneous boundary value problems
(BVPs) of the form
\begin{equation}
\label{eq:basic}
\begin{array}{r ll}
-\Delta u (\vct{x}) + b(\vct{x})u_x(\vct{x})+c(\vct{x}) u_y(\vct{x}) +d(\vct{x}) u(\vct{x})&\!\!\!\! =0 \quad &\vct{x} \in \Omega \\
u (\vct{x})  &\!\!\!\! = g(\vct{x}) \quad &\vct{x} \in \Gamma,
\end{array}
\end{equation}
where $\Omega = [0,1]^{2}$ is the unit square in $\mathbb{R}^{2}$,
where $\Gamma$ is the boundary of $\Omega$,
and where $b$, $c$, and $d$ are functions on $\Omega$.
We assume that the only information sought is the normal derivative
of $u$ at $\Gamma$. In other words, the objective is to construct an approximation
to the Dirichlet-to-Neumann operator associated with the elliptic differential
operator in (\ref{eq:basic}).


The proposed solver is particularly efficient for situations where (\ref{eq:basic})
needs to be solved for a sequence of different boundary data functions $g$.
The solver has two steps: (1) Build the approximate Dirichlet-to-Neumann operator
for a given set of functions $b$, $c$, and $d$. (2) Determine the Neumann data
$\partial u/\partial n$ for any given Dirichlet data $g$
by applying the approximate Dirichlet-to-Neumann operator.
The key claim of the paper is that the ``build stage'' can be executed in $O(N)$
operations, and the ``solve stage'' can be executed in $O(N^{0.5})$ operations.
In contrast, classical nested dissection requires $O(N^{1.5})$ and $O(N)$
operations for the two steps, respectively.

\subsection{Motivation}
While the present paper addresses the specific BVP (\ref{eq:basic}) on a simple square
domain, the technique can be extended for building solution operators to elliptic
boundary value problems of the form
\begin{equation}
\label{eq:BVP_general}
\left\{\begin{array}{lll}
Au(\vct{x}) =&\!\!\!\!\!\! f(\vct{x})\qquad&\vct{x} \in \Omega\\
Bu(\vct{x}) =&\!\!\!\!\!\! g(\vct{x})\qquad&\vct{x} \in \Gamma,
\end{array}\right.
\end{equation}
where $\Omega$ is a domain in $\mathbb{R}^{2}$ or $\mathbb{R}^{3}$ with boundary
$\Gamma$, where $A$ is an elliptic partial differential operator, and where $B$
is a trace operator (representing boundary conditions like Dirichlet, Neumann, mixed, etc.).
The goal of this paper is to illustrate that when the body load $f$ is zero
and the solution $u$ and/or its derivatives are sought only near the boundary, the relevant solution operator
can be constructed at moderate cost, and applied almost instantaneously. This
opens up the possibility of high accuracy computational simulations to be carried
out in real time for 3D problems such as elasticity involving composite materials,
electrostatics in domains with variable conductivity, acoustic and electromagnetic
scattering problems (at long and intermediate wave-lengths at least), and many others.

In some applications such as seismic testing and automatic multilevel substructuring (AMLS),
there is a small number of localized body loads inside of the domain.  This paper
details how the solution operator can be found with an increased cost which is still less
than classic techniques.

\subsection{Discretization}
The method described is applicable
to a variety of geometries and discretization schemes (finite elements, finite
differences, etc.). For simplicity of presentation, we restrict our attention
to the model problem where a square domain $\Omega$ is discretized via a finite
difference scheme on a regular $n\times n$ square mesh. The resulting linear
system takes the form
\begin{equation}
\mtx{A} \vct{u} = \vct{b}
\label{eq:basic_nd}
\end{equation}
where $\mtx{A}$ is an $n^{2} \times n^{2}$ sparse matrix.
We let $N = n^{2}$ denote the total number of grid points.

\vspace{2mm}

\noindent\textit{\small\textbf{Example:}
When (\ref{eq:basic}) represents the Laplace equation ($b = c = d = 0$)
and the standard five-point finite difference stencil is used in the
discretization, $\mtx{A}$
consists of $n\times n$ blocks, each of size $n\times n$,
\begin{equation*}
\mtx{A} = \left[\begin{array}{rrrrr}
           \mtx{B} & -\mtx{I} & 0 & 0 &\cdots \\
          -\mtx{I} & \mtx{B} & -\mtx{I} & 0 & \cdots\\
          0&-\mtx{I} & \mtx{B} & -\mtx{I} & \cdots \\
         0& 0&-\mtx{I} & \mtx{B} & \cdots  \\
\vdots & \vdots & \vdots & \vdots &
           \end{array}\right],
\quad\mbox{where}\quad
\mtx{B} = \left[\begin{array}{rrrrr}
           4 & -1 & 0 & 0 &\cdots \\
          -1 & 4 & -1 & 0 & \cdots\\
          0&-1 & 4 & -1 & \cdots \\
         0& 0&-1 & 4 & \cdots  \\
\vdots & \vdots & \vdots & \vdots &
           \end{array}\right],
\end{equation*}
and where $\mtx{I}$ is the $n\times n$ identity matrix.}

\subsection{Existing fast solvers}
There already exist many efficient techniques for solving (\ref{eq:basic_nd}),
including:

\vspace{1mm}

\noindent
\textbf{\textit{Iterative methods:}} These techniques construct a
sequence of successively more accurate approximate solutions by applying
the matrix $\mtx{A}$ to a sequence of vectors.
Since the $N\times N$ matrix $\mtx{A}$ has $O(N)$ non-zero entries,
the resulting solver has $O(N)$ complexity whenever convergence is fast.
It is difficult to predict the convergence rate of iterative methods and often
a customized pre-conditioner is required to accelerate the schemes.

\vspace{1mm}

\noindent
\textbf{\textit{Multigrid methods:}} These techniques can be viewed as a special case of iterative
methods. They can in certain circumstances reach very high performance by
decomposing the matrix in a sequence of different scales; since the matrix
is well-conditioned on each scale, very fast convergence often results.

\vspace{1mm}

\noindent
\textbf{\textit{Direct methods:}} Direct solvers (such as Gaussian elimination) which compute
a solution in a single shot are considered more stable and robust than iterative methods.
Proper ordering of the nodes often allows Gaussian elimination to be
executed at $O(N^{1.5})$ complexity (\cite{george_1973}), and the resulting ``nested
dissection'' approach is quite competitive for moderate problem sizes
(up to about $N \sim 10^{6}$). More recently, it has been shown that by
exploiting additional structure in the coefficient matrix, the nested
dissection method can be accelerated to (close to) linear complexity, see, e.g.,
\cite{2009_xia_multifrontal,2009_grasedyck_FEM,2010_ying_nesteddissection}.


\subsection{Novelty of the present work and comparison to existing methodology}
The proposed solver is based on the classical nested dissection algorithm of
\cite{george_1973}. The key distinction to classical nested dissection
is that special structure in the dense so called ``frontal matrices'' are
exploited to reduce the cost of the pre-computation from $O(N^{1.5})$ to $O(N)$,
and the cost of the solve from $O(N)$ to $O(N^{0.5})$. To be precise, we
approximate off-diagonal blocks of the dense frontal matrices by low-rank
matrices; we do this using the structured matrix format described in
\cite{2010_jianlin_fast_hss,m2011_1D_survey}, which can be viewed as a variation
of the well established ``$\mathcal{H}$-matrix'' and ``$\mathcal{H}^{2}$-matrix''
formats of Hackbusch and co-workers (see, e.g., \cite{hackbusch,2007_borm_H2b}).

The observation that the dense matrix computations involving the frontal matrices
can be accelerated using structured matrix algebra has recently been made
in, e.g., \cite{2009_xia_multifrontal,2009_grasedyck_FEM,m2009_onion,2010_ying_nesteddissection}.
Our work is slightly different in that it is based on the hierarchical construction
of Schur complements, and directly leads to a discrete approximation of the
Dirichlet-to-Neumann operator on the full domain. This greatly simplifies the
construction of the boundary-to-boundary solution operator. It also leads to
algorithms that can readily handle a problem involving a sparse body load
(see Section \ref{sec:bodyload}).

While the present work considers only a regular square grid in 2D, the method
can be extended to more general grids, and to 3D problems. A discussion of the
expected performance in these cases can be found in Section \ref{sec:conc}.
Our key claims regarding asymptotic complexity are summarized in Table \ref{tab:scaling}.

An early version of the work reported appeared in the Ph.D.~dissertation \cite{Adiss}.

\begin{remark}
State-of-the-art iterative solvers such as, e.g., multigrid will sometimes outperform
the new accelerated nested dissection technique for a stand-alone solve.
However, the new solver is much faster for subsequent solves;
its asymptotic cost is only $O(N^{0.5})$ and its practical efficiency is such
that a problem on a grid with $4000 \times 4000$ nodes can be solved in only
0.1 seconds on a standard office laptop. Moreover, since the solver is direct,
it handles with ease many problems that are challenging for iterative methods
(including multigrid), such as for instance vibration problems in situations
where the domain is much larger than the wave-length.
\end{remark}

\begin{table}[ht]
\centering\small
\begin{tabular}{|c|ll|ll|ll|}
\hline
   & \multicolumn{2}{l|}{Build solution operators} &
     \multicolumn{2}{l|}{Solve with no body load}  &
     \multicolumn{2}{l|}{Solve with general body load} \\ \hline
2D & $N$\qquad & ($N^{3/2}$)     &  $N^{1/2}$  & ($N$)       & $N$              & ($N \log N$)      \\ \hline
3D & $N^{4/3}$ & ($N^{2}$)       &  $N$ & ($N^{4/3}$) & $N\log N$  & ($N^{4/3} \log N$) \\ \hline
\end{tabular}
\caption{\label{tab:scaling}
Summary of asymptotic costs of the proposed direct solver in two and three dimensions.
For comparison, the costs of classical nested dissection are given in parenthesis.
For special (e.g.~sparsely supported) body loads, better asymptotics than those listed in the table can be achieved,
see Section \ref{sec:bodyload}.
The asymptotics given for 3D problems are predictions for a simplistic generalization of the proposed scheme;
it is likely that these numbers could be further improved.}
\end{table}

\subsection{Outline of paper}
The paper describes an $O(N)$ variation of the nested dissection method which
computes the global Dirichlet-to-Neumann operator.
Section \ref{sec:quadtree} describes a hierarchical partitioning of the grid into
a quad-tree of nested boxes.
Section \ref{sec:merge} describes a variation of the classical nested dissection
technique that computes a hierarchy of solution operators for each box in the quad-tree.
The solution operators have internal structure (see Section \ref{sec:HBS}) which is
exploited to improve the complexity of the method from $O(N^{1.5})$ to $O(N)$, see
Sections \ref{sec:fastHBS} and \ref{sec:accel}.
Section \ref{sec:num} reports the results of numerical experiments that substantiate
our claims on the asymptotic complexity and accuracy of the method.

\section{Tree structure}
\label{sec:quadtree}

The direct solver described in this note is based on the classical nested dissection algorithm,
and uses an analogous (but not identical) tree structure on the computational grid. This section
formally defines the tree structure for our simple model geometry.

Let $\Omega$ denote the square domain introduced in Section \ref{sec:intro}, and suppose
that it is discretized using a uniform $n\times n$ grid. Let $N = n^{2}$ denote the number
of points in the grid, and let $N_{\rm leaf}$ denote a tuning parameter
chosen so that a matrix of size $N_{\rm leaf}\times N_{\rm leaf}$ can be inverted quickly by brute force.
The optimal choice of $N_{\rm leaf}$ depends on the computing environment,
but we have found that $N_{\rm leaf} = 4096$ is often a good choice.
Let $L$ be the smallest integer such that when $\Omega$ is partitioned
into $4^L$ equisized boxes, each box contains no more than $N_{\rm leaf}$ points.
These $4^{L}$ small boxes are called the \textit{leaves} of the tree.
Merge the leaves by sets of fours into boxes with twice the side length,
to form the $4^{L-1}$ boxes that make up the next level in the tree.
This process is repeated until $\Omega$ is recovered.
We call $\Omega$ the \textit{root} of the tree.

The set consisting of all boxes of the same size forms what we
call a \textit{level}. We label the levels using the integer $\ell = 0,\,1,\,2,\,\dots,\,L$,
with $\ell = 0$ denoting the root, and $\ell = L$ denoting the leaves.


\section{A variation of the nested dissection algorithm}
\label{sec:merge}

This section describes a direct solver
that is particularly fast
for what we call ``pure'' boundary value problems such as (\ref{eq:basic}) in which
there is no body load, and where the solution is sought only near the boundary.
The idea is to construct a solution operator $\mtx{G}$ that maps the given boundary
data to the sought potential values (or flows) on the boundary. Letting $N_{\rm b}$
denote the number of nodes on the boundary of the domain, $\mtx{G}$ is a dense
$N_{\rm b} \times N_{\rm b}$ matrix.

Technically, the solution operator $\mtx{G}$ is constructed via a divide-and-conquer
approach (analogous to the one used in the classical nested dissection scheme): First
a solution operator is constructed for each ``leaf'' in the quadtree described in
Section \ref{sec:quadtree}, then solution operators for larger boxes are constructed
via a hierarchical merging process in a single sweep through the tree, going from
smaller to larger boxes.

For a grid with $N$ nodes, the process described in this section requires
$O(N^{1.5})$ operations to construct the solution operator, and then each
subsequent solve (which consists merely of applying the solution operator)
requires $O(N)$ operations. Techniques for accelerating these two costs to
$O(N)$ and $O(N^{0.5})$, respectively, are then described in Sections \ref{sec:HBS}, \ref{sec:fastHBS}
and \ref{sec:accel}.

\subsection{The solution operator and the Schur complement}
\label{sec:schur}
\begin{figure}[ht]
\begin{center}
\begin{tabular}{ccc}
\setlength{\unitlength}{1mm}
\begin{picture}(45,40)
\put(00,00){\includegraphics[width=45mm]{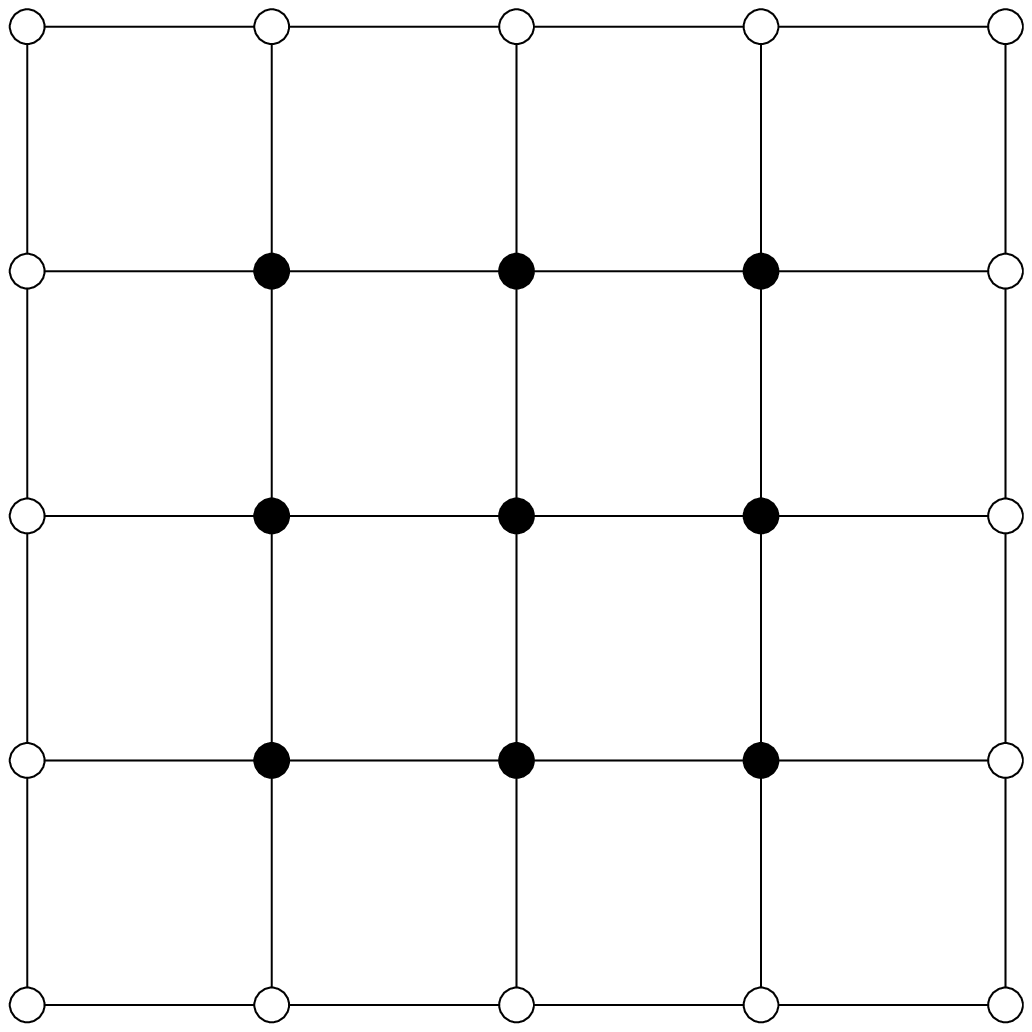}}
\put(-4,14){$P_{\rm b}$}
\put(27,15){$P_{\rm i}$}
\end{picture}
&&
\setlength{\unitlength}{1mm}
\begin{picture}(45,40)
\put(00,00){\includegraphics[width=45mm]{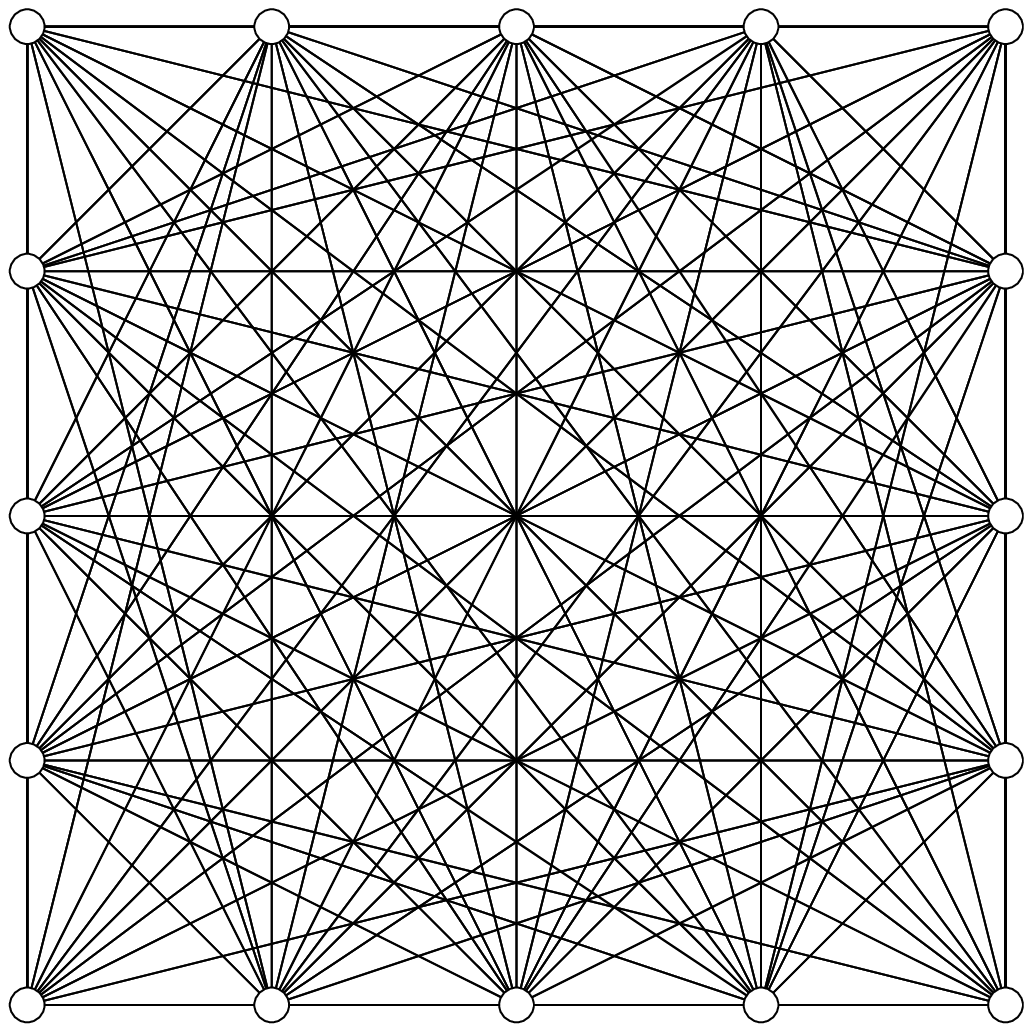}}
\end{picture}\\
(a) &\mbox{}\hspace{20mm}\mbox{}& (b)
\end{tabular}
\end{center}
\caption{(a) Labeling of nodes for constructing the Schur complement
of a leaf. $P_{\rm i}$ are the interior nodes
(solid), and $P_{\rm b}$ are the boundary nodes (hollow). (b) After the merge, all internal
nodes are ``eliminated'' but now all nodes communicate directly (i.e.~the Schur complement
$\mtx{S}$ is dense).}
\label{fig:onebox}
\end{figure}

This subsection provides a precise definition of the concept of a
``solution operator'' associated
with a subdomain $P$ of the computational grid. For simplicity, we assume that $P$
is a square or rectangular domain. We partition $P$ into interior nodes and boundary nodes:
$$
P = P_{\rm i} \cup P_{\rm b},
$$
where $P_{\rm i}$ is defined as the set of nodes that have all four neighbors
inside $P$, see Figure \ref{fig:onebox}. (Note that the set $P$ consists of all nodes
at which the potential is unknown, and $P_{\rm b}$ is the outermost ring of these nodes,
\textit{not} the nodes at which Dirichlet data is prescribed.)

Let $\vct{u}_{\rm b}$ and $\vct{u}_{\rm i}$ denote the potentials at the boundary
nodes and the interior nodes, respectively. Reordering the equilibrium equation
(restricted to $P$), we find that $\vct{u}_{\rm b}$ and $\vct{u}_{\rm i}$ must satisfy
\begin{equation}
\label{eq:chalk2}
\left[\begin{array}{cc}
\mtx{A}_{\rm b,b} & \mtx{A}_{\rm b,i} \\
\mtx{A}_{\rm i,b} & \mtx{A}_{\rm i,i}
\end{array}\right]\,
\left[\begin{array}{c}
\vct{u}_{\rm b} \\
\vct{u}_{\rm i}
\end{array}\right] =
\left[\begin{array}{c}
\vct{f}_{\rm b} \\
\vct{f}_{\rm i}
\end{array}\right].
\end{equation}
Eliminating $\vct{u}_{\rm i}$ from (\ref{eq:chalk2}), we find
$$
\mtx{S}\,\vct{u}_{\rm b} = \vct{f}_{\rm b}- \mtx{A}_{\rm b,i}\mtx{A}^{-1}_{\rm i,i}\vct{f}_{\rm i},
$$
where $\mtx{S}$ is the matrix
\begin{equation}
\label{eq:def_schur_nd}
\mtx{S} = \mtx{A}_{\rm b,b} - \mtx{A}_{\rm b,i}\,\mtx{A}_{\rm i,i}^{-1}\,\mtx{A}_{\rm i,b}.
\end{equation}
We refer to $\mtx{S}$ as the \textit{Schur complement} associated with the subdomain $P$;
the \textit{solution operator} is then $\mtx{G} = \mtx{S}^{-1}$.
In the case of no body loads $\vct{f}_{\rm i} = \vct{0}$, thus the update to the
right hand side on the boundary is not necessary.

\subsection{Merging two Schur complements}
\label{sec:mergetwo}
\begin{figure}[ht]
\begin{center}
\begin{tabular}{ccc}
\setlength{\unitlength}{1mm}
\begin{picture}(47,24)
\put(00,00){\includegraphics[height=22mm]{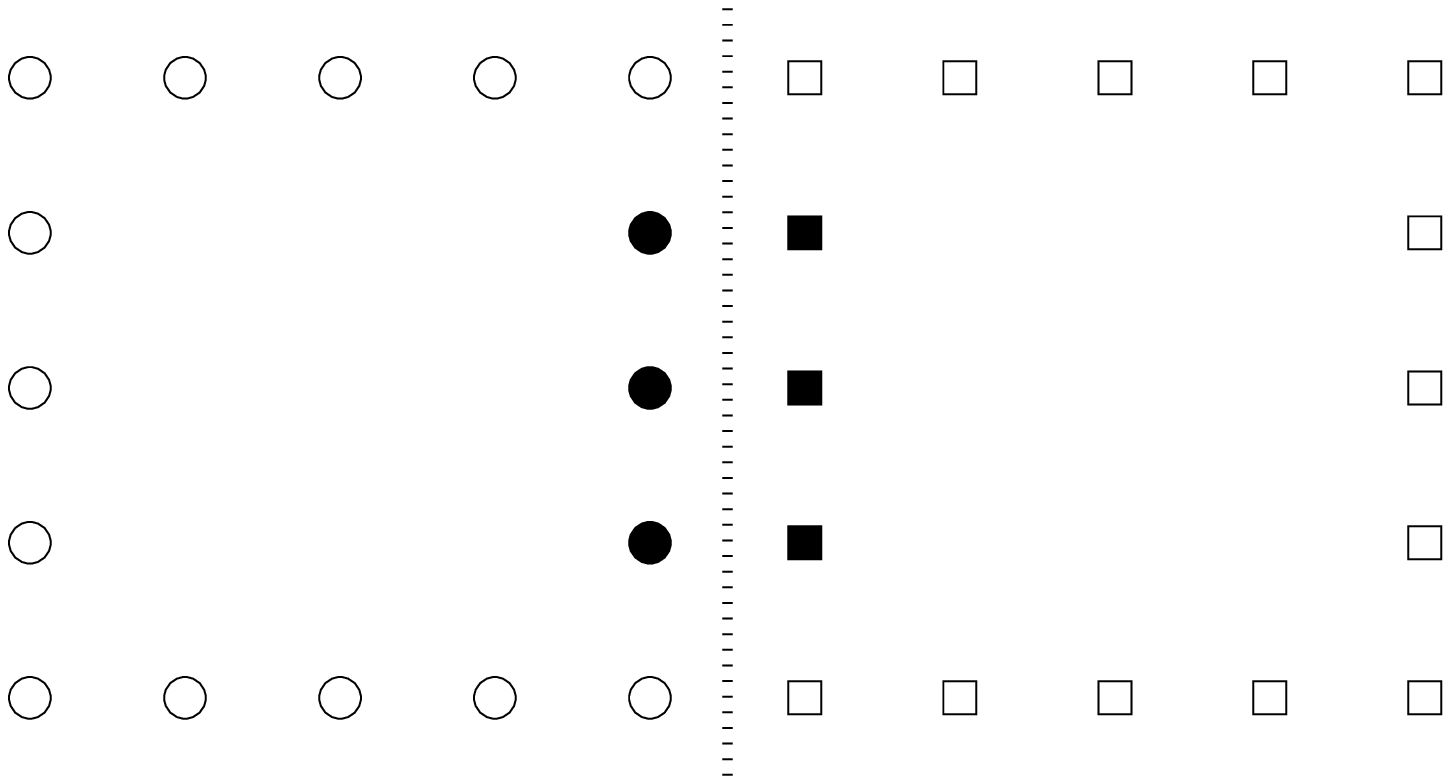}}
\put(-5,11){$P_{1}$}
\put(13,11){$P_{3}$}
\put(25,11){$P_{4}$}
\put(41,11){$P_{2}$}
\put(08,23){$\Omega^{\rm(w)}$}
\put(28,23){$\Omega^{\rm(e)}$}
\end{picture}
&
\setlength{\unitlength}{1mm}
\begin{picture}(47,24)
\put(00,00){\includegraphics[height=22mm]{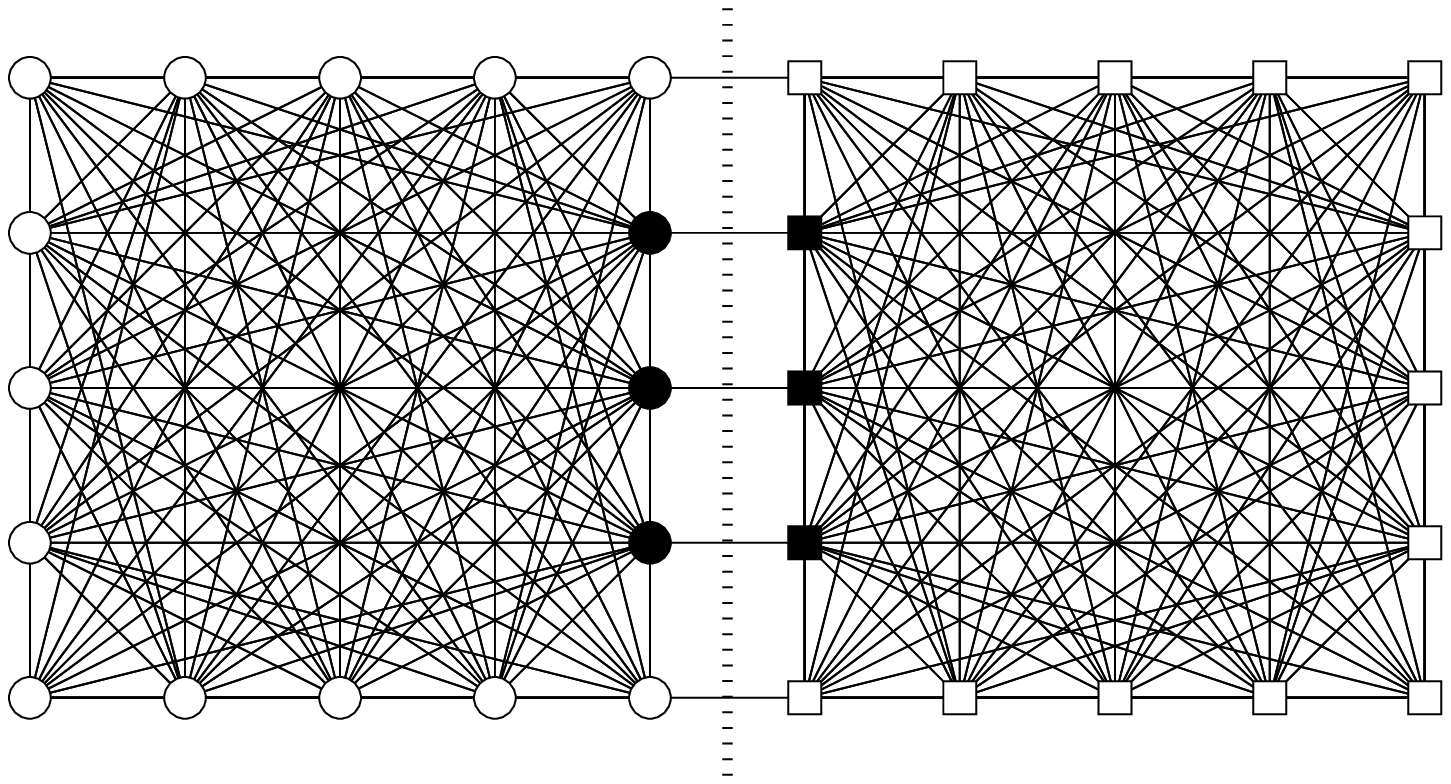}}
\end{picture}
&
\setlength{\unitlength}{1mm}
\begin{picture}(47,24)
\put(00,00){\includegraphics[height=22mm]{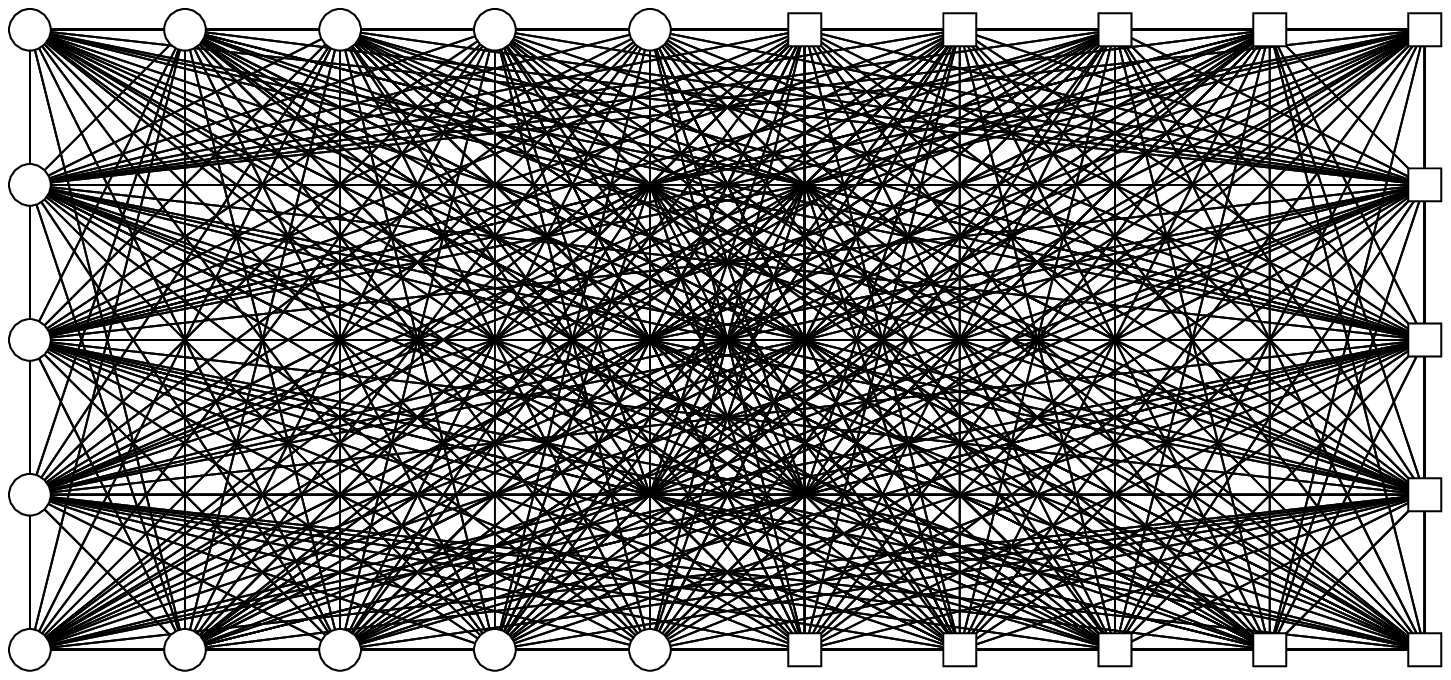}}
\end{picture}
\\
(a) & (b) & (c)
\end{tabular}
\end{center}
\caption{(a) Labeling of nodes for the merge operation described in Section \ref{sec:merge}.
The nodes in $P_{1}$ and $P_{3}$ are round, and the nodes in $P_{2}$ and $P_{4}$ are square.
The solid nodes are interior to the union of the two boxes $\Omega^{\rm(w)}$ and $\Omega^{\rm(e)}$.
(b) Connections between nodes before the merge.
(c) Connections between nodes after eliminating the interior (solid) nodes.
}
\label{fig:mergetwo}
\end{figure}

In this section, we present a technique for merging the Schur complements
for two adjacent boxes. Let us call the two boxes $\Omega^{\rm(w)}$ and
$\Omega^{\rm(e)}$ (for ``west'' and ``east'').
Further, let $P^{\rm(w)}$
and $P^{\rm(e)}$ denote the nodes on the boundaries of these two boxes,
and let $\mtx{S}^{\rm(w)}$ and $\mtx{S}^{\rm(e)}$ denote Schur complements
supported on these two sets of boundary nodes, see Figure \ref{fig:mergetwo} (a).

The objective of the merge is to eliminate the nodes that are now
``interior'' to the larger box formed by the union of the two smaller boxes;
these nodes are marked as blue in Figure \ref{fig:mergetwo} (b).
To eliminate these points, we first partition the nodes in $P^{\rm(w)}$ and
$P^{\rm(e)}$ so that
\begin{equation}
\label{eq:partition_we}
P^{\rm(w)} = P_{1} \cup P_{3},
\qquad\mbox{and}\qquad
P^{\rm(e)} = P_{2} \cup P_{4},
\end{equation}
and so that $P_{1} \cup P_{2}$ forms the boundary of the larger box,
while the nodes in $P_{3}\cup P_{4}$ are interior, see Figure \ref{fig:mergetwo} (c).
Partition the Schur complements $\mtx{S}^{\rm(w)}$ and $\mtx{S}^{\rm(e)}$ analogously:
$$
\mtx{S}^{\rm(w)} = \left[\begin{array}{cc}
\mtx{S}_{11} & \mtx{S}_{13} \\
\mtx{S}_{31} & \mtx{S}_{33}
\end{array}\right],
\qquad\mbox{and}\qquad
\mtx{S}^{\rm(e)} = \left[\begin{array}{cc}
\mtx{S}_{22} & \mtx{S}_{24} \\
\mtx{S}_{42} & \mtx{S}_{44}
\end{array}\right].
$$

Supposing that the right hand side has been updated to account for any interior body loads, equation (\ref{eq:basic_nd}) restricted
to the union of the two boxes now reads

\begin{equation}
\label{eq:mergetwo_nd}
\left[\begin{array}{cc|cc}
\mtx{S}_{11} & \mtx{A}_{12} & \mtx{S}_{13} & \mtx{A}_{14}      \\
\mtx{A}_{21} & \mtx{S}_{22} & \mtx{A}_{23}      & \mtx{S}_{24} \\ \hline
\mtx{S}_{31} & \mtx{A}_{32}     & \mtx{S}_{33} & \mtx{A}_{34} \\
\mtx{A}_{41}   & \mtx{S}_{24} & \mtx{A}_{43} & \mtx{S}_{44}
\end{array}\right]\,
\left[\begin{array}{c}
\vct{u}_{1} \\ \vct{u}_{2} \\ \vct{u}_{3} \\ \vct{u}_{4}
\end{array}\right] =
\left[\begin{array}{c}
\vct{f}_{1} \\ \vct{f}_{2} \\ \vct{f}_{3} \\ \vct{f}_{4}
\end{array}\right],
\end{equation}
where $\mtx{A}_{ij}$ are the relevant sub-matrices of the original discrete Laplacian $\mtx{A}$.
From (\ref{eq:mergetwo_nd}), one finds that the Schur complement of the union box is
\begin{equation}
\label{eq:S_fromtwo_nd}
\mtx{S} =
\left[\begin{array}{cc}
\mtx{S}_{11} & \mtx{A}_{12} \\
\mtx{A}_{21} & \mtx{S}_{22}
\end{array}\right] -
\left[\begin{array}{cc}
\mtx{S}_{13} & \mtx{A}_{14}    \\
\mtx{A}_{23}  & \mtx{S}_{24}
\end{array}\right]\,
\left[\begin{array}{cc}
\mtx{S}_{33} & \mtx{A}_{34} \\
\mtx{A}_{43} & \mtx{S}_{44}
\end{array}\right]^{-1}\,
\left[\begin{array}{cc}
\mtx{S}_{31} & \mtx{A} _{32}     \\
\mtx{A}_{41}      & \mtx{S}_{42}
\end{array}\right].
\end{equation}
The updated right hand side is
\begin{equation*}
\vtwo{\vct{f}_{1}}{\vct{f}_{2}} -
\left[\begin{array}{cc}
\mtx{S}_{13} & \mtx{A}_{14}    \\
\mtx{A}_{23}  & \mtx{S}_{24}
\end{array}\right]\,
\left[\begin{array}{cc}
\mtx{S}_{33} & \mtx{A}_{34} \\
\mtx{A}_{43} & \mtx{S}_{44}
\end{array}\right]^{-1}\,
\vtwo{\vct{f}_{3}}{\vct{f}_{4}}.\end{equation*}

\begin{remark}
The matrices $\mtx{A}_{14}$, $\mtx{A}_{23}$, $\mtx{A}_{32}$, and $\mtx{A}_{41}$ are typically
very sparse. In fact, when equation (\ref{eq:basic}) is discretized with a $5$-point stencil,
these matrices are identically zero.
\end{remark}

\begin{figure}[ht]
\begin{tabular}{cc}
 \setlength{\unitlength}{1mm}
 \begin{picture}(55,55)
 \put(-15,00){\includegraphics[width=80mm]{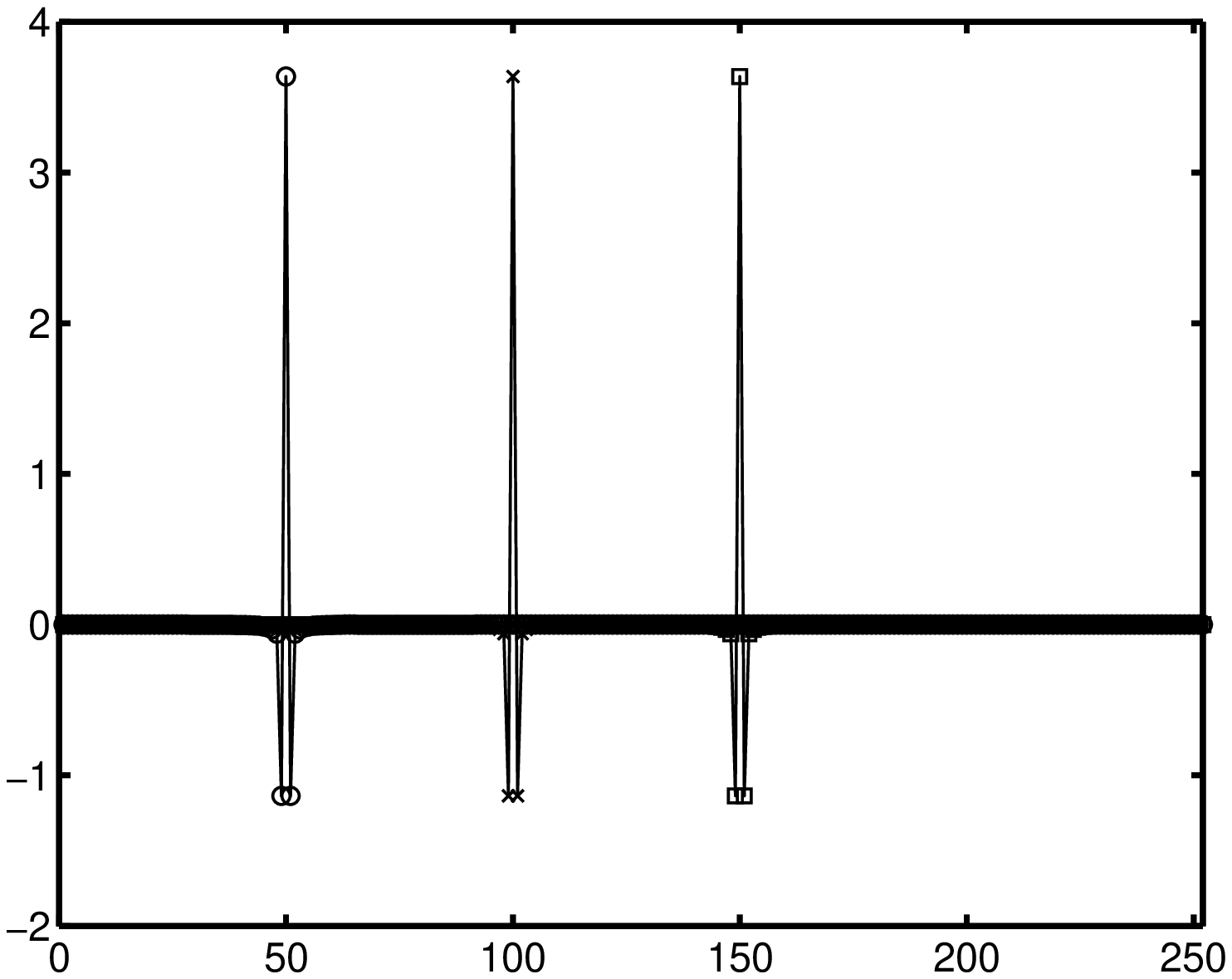}}
 \put(25,-1){(a)}
 \end{picture}
&
 \setlength{\unitlength}{1mm}
 \begin{picture}(55,55)
 \put(00,00){\includegraphics[width=80mm]{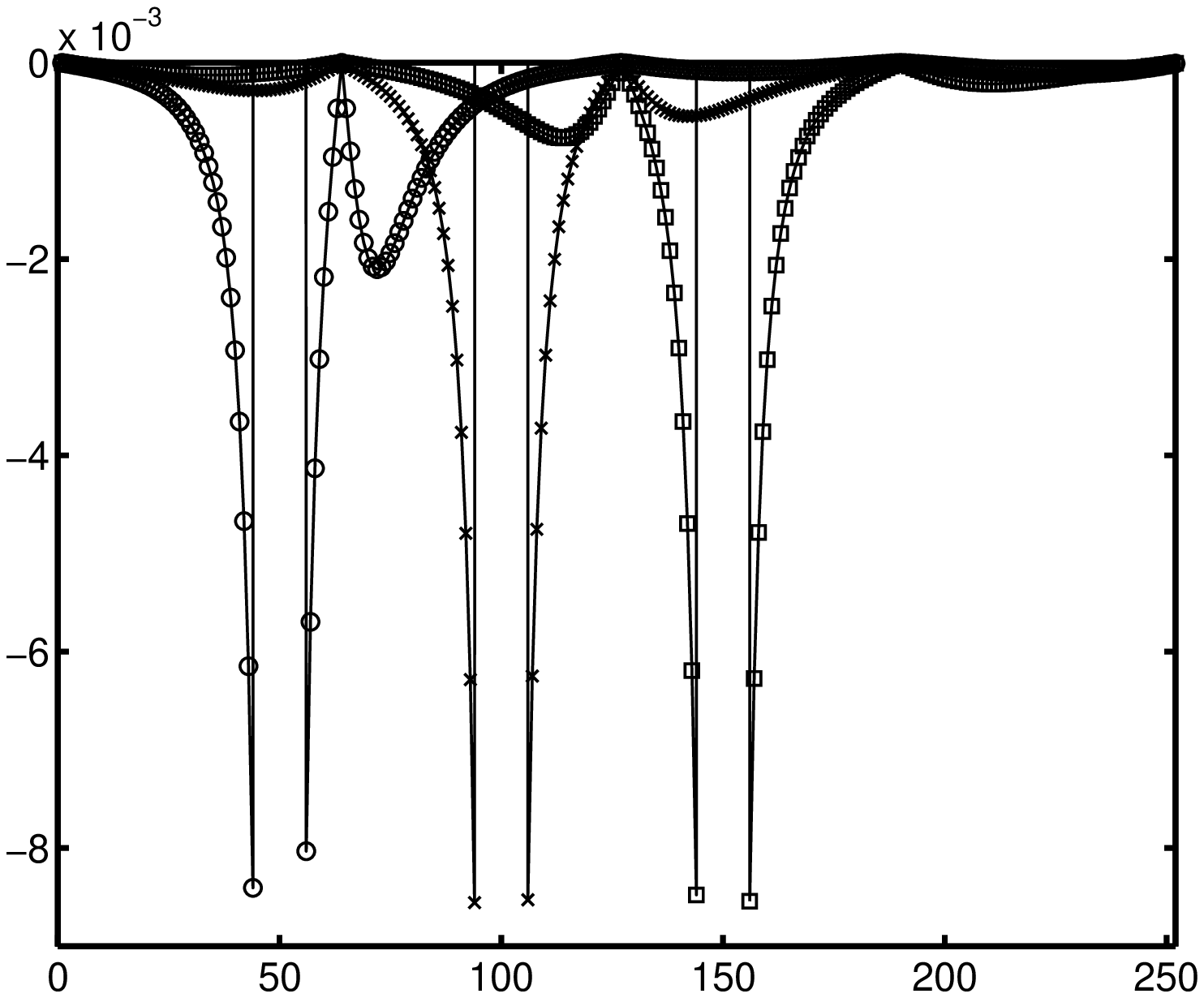}}
 \put(40,-1){(b)}
 \end{picture}
\\
\end{tabular}
\caption{(a) Three rows of a typical Schur complement for a Laplace problem with $252$ points on the boundary.
The $50^{\rm th}$ ($\circ$), $100^{\rm th}$ ($\times$), and $150^{\rm th}$ ($\Box$) rows are shown. Note how the matrix
is completely dominated by the elements close to the diagonal.
(b) The same plot but with a different scale on the vertical axis to show the smoothness in the far-field.}
\label{fig:row_plot}
\end{figure}

\begin{remark}
While the Schur complement is a dense matrix (cf.~Figure \ref{fig:mergetwo}),
the interactions between distant points can to high precision be approximated by low rank
matrices. This property can be conjectured by inspecting a computed Schur complement and
observing that each row is smooth away from the diagonal. Figure (\ref{fig:row_plot}) illustrates this point.
The figure also illustrates that the Schur complement is strongly diagonal dominant,
which is consistent with the fact that it is a discrete analog of the Dirichlet-to-Neumann
operator, which is a hyper-singular integral operator (it \textit{reduces} the smoothness
of any boundary function it operates on in a manner similar to a differentiation operator).
\end{remark}

\subsection{The full algorithm}
\label{sec:full}
For future reference, let us summarize the algorithm described:

\vspace{1.5mm}

\begin{enumerate}
\item \textit{Construct a quad-tree:} Partition the grid into a hierarchy of boxes
as described in Section \ref{sec:quadtree}.

\vspace{1.5mm}

\item \textit{Process the leaves:} For each leaf box in the tree, construct its
Schur complement as described in Section \ref{sec:schur}.  If the box has body
loads, update the right hand side.

\vspace{1.5mm}

\item \textit{Hierarchical merge:} Loop over all levels of the tree, from finer to
smaller. For each box on a level, compute its Schur complement by merging the (already
computed) Schur complements of its children.
Note that the merge of four children can be executed via three of the
pair-wise merges described in Section \ref{sec:mergetwo}.
If the interior points have body loads, update the right hand side.

\vspace{1.5mm}

\item \textit{Process the root of the tree:} After completing Step 3, the Schur
complement for the entire domain is available. Invert (or factor) it to construct
the solution operator.
\end{enumerate}

\begin{remark}
For simplicity, the algorithm is described in a level-by-level manner
(process all leaves first, then proceed one level at a time in going upwards).
In fact, there is flexibility to travel through the tree in any order that ensures
that no node is processed before its children. Since all Schur complements can be
discarded once their information has been passed on to a parent, smarter orderings
can greatly reduce the memory requirements \cite{ess76-element}.
\end{remark}

\subsection{Asymptotic complexity of the algorithm}
As before, let $N = n^2$ denote the total number of points in the grid, let $N_{\rm leaf}$
denote the maximum number of points on a leaf, and let $L$ denote the number of levels
so that $N \leq 4^{L}\,N_{\rm leaf}$.

The cost to process one leaf in Step 2 in Section \ref{sec:full} is $O(N_{\rm leaf}^{2})$
(exploiting that the matrix $\mtx{A}_{\rm i,i}$ in (\ref{eq:def_schur_nd}) is band-diagonal).
Since there are $4^{L}$ leaves, the total cost of Step 2 is therefore
$4^{L}\,N_{\rm leaf}^{2} \sim N\,N_{\rm leaf}$. Since $N_{\rm leaf}$ is a small
constant number (in principle one could set $N_{\rm leaf}=1$) the leaf processing cost is $O(N)$.

Next consider the cost of constructing the Schur complement of a box on level $\ell$ in
executing Step 3 in Section \ref{sec:full}. Note that all boxes involved have $O(n\,2^{-\ell})$
points along each side. Since some matrices in (\ref{eq:S_fromtwo_nd}) are dense, the
cost for each merge is proportional to $(n\,2^{-\ell})^{3} = n^{3}\,2^{-3\ell}$.
Since we need to compute $2^{2\ell}$ Schur complements on level $\ell$, the total cost
of Step 3 is then $\sum_{\ell = 1}^{L}2^{2\ell}\,n^{3}\,2^{-3\ell} =
n^{3}\sum_{\ell = 1}^{L}2^{-\ell} \sim n^{3}$.

Since the cost of the final inversion/factorization in Step 4 is $O(n^{3})$,
the total cost of the algorithm in Section \ref{sec:full} is $O(n^{3}) = O(N^{1.5})$.

\section{Compressible matrices}
\label{sec:HBS}
To improve the scaling of the nested dissection method, a more efficient technique
for evaluating (\ref{eq:S_fromtwo_nd}) will be implemented. We will exploit that
while the matrices $\mtx{S}_{ij}$ are all dense, they in the present context have
additional structure: $\mtx{S}_{ij}$ is when $i \neq j$ to high precision rank deficient,
and $\mtx{S}_{ii}$ has a structure that we call \textit{Hierarchically Block Separable (HBS)}.
This section briefly describes the HBS property, for details see \cite{m2011_1D_survey}.
We note that the HBS property is very similar to the concept of \textit{Hierarchically
Semi-Separable (HSS)} matrices \cite{2007_shiv_sheng,gu_divide} which has previously
been used in an analogous context \cite{2009_xia_multifrontal}. Other researchers have
used the somewhat related $\mathcal{H}$-matrix concept for similar purposes
\cite{2009_grasedyck_FEM,2010_ying_nesteddissection}.

\subsection{Block separable}
Let $\mtx{H}$ be an $mp\times mp$ matrix that is blocked into $p\times p$ blocks,
each of size $m\times m$.

We say that $\mtx{H}$ is ``block separable'' with ``block-rank'' $k$
if for $\tau = 1,\,2,\,\dots,\,p$, there exist $m\times k$
matrices $\mtx{U}_{\tau}$ and $\mtx{V}_{\tau}$ such that each off-diagonal
block $\mtx{H}_{\sigma,\tau}$ of $\mtx{H}$ admits the factorization
\begin{equation}
\label{eq:yy1}
\begin{array}{cccccccc}
\mtx{H}_{\sigma,\tau}  & = & \mtx{U}_{\sigma}   & \tilde{\mtx{H}}_{\sigma,\tau}  & \mtx{V}_{\tau}^{*}, &
\quad \sigma,\tau \in \{1,\,2,\,\dots,\,p\},\quad \sigma \neq \tau.\\
m\times m &   & m\times k & k \times k & k\times m
\end{array}
\end{equation}

Observe that the columns of $\mtx{U}_{\sigma}$ must form a basis for
the columns of all off-diagonal blocks in row $\sigma$, and
analogously, the columns of $\mtx{V}_{\tau}$ must form a basis for the
rows in all the off-diagonal blocks in column $\tau$. When (\ref{eq:yy1})
holds, the matrix $\mtx{H}$ admits a block factorization
\begin{equation}
\label{eq:yy2}
\begin{array}{cccccccccc}
 \mtx{H} &  =& \mtx{U}&\tilde{\mtx{H}}& \mtx{V}^{*} & +&  \mtx{D},\\
mp\times mp &   & mp\times kp & kp \times kp & kp\times mp && mp \times mp\\
\end{array}
\end{equation}
where
$$
\mtx{U} = \mbox{diag}(\mtx{U}_{1},\,\mtx{U}_{2},\,\dots,\,\mtx{U}_{p}),\quad
\mtx{V} = \mbox{diag}(\mtx{V}_{1},\,\mtx{V}_{2},\,\dots,\,\mtx{V}_{p}),\quad
\mtx{D} = \mbox{diag}(\mtx{D}_{1},\,\mtx{D}_{2},\,\dots,\,\mtx{D}_{p}),
$$
and
$$\tilde{\mtx{H}} = \left[\begin{array}{cccc}
0 & \tilde{\mtx{H}}_{12} & \tilde{\mtx{H}}_{13} & \cdots \\
\tilde{\mtx{H}}_{21} & 0 & \tilde{\mtx{H}}_{23} & \cdots \\
\tilde{\mtx{H}}_{31} & \tilde{\mtx{H}}_{32} & 0 & \cdots \\
\vdots & \vdots & \vdots
\end{array}\right].
$$

\subsection{Heirarchically Block-Separable}
Informally speaking, a matrix $\mtx{H}$ is \textit{Heirarchically Block-Separable} (HBS),
if it is amenable to a \textit{telescoping} block factorization.  In other words,
in addition to the matrix $\mtx{H}$ being block separable, so is $\tilde{\mtx{H}}$
once it has been reblocked to form a matrix with $p/2 \times p/2$ blocks.
Likewise, the middle matrix from the block separable factorization
of $\tilde{\mtx{H}}$ will be block separable, etc.

In this section, we describe properties and the factored representation of HBS matrices.
Details on constructing the factorization are provided in \cite{m2011_1D_survey}.

\subsubsection{A binary tree structure}
\label{sec:tree}
\begin{figure}
\setlength{\unitlength}{1mm}
\begin{picture}(169,41)
\put(20, 0){\includegraphics[height=41mm]{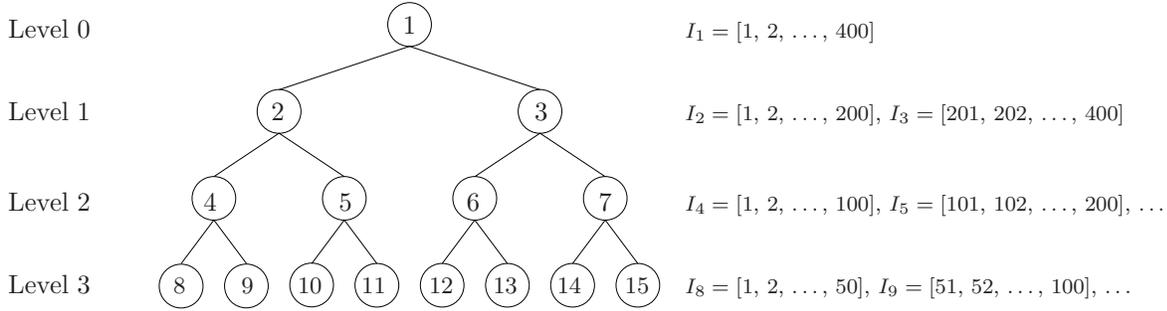}}
\put( 0,36){Level $0$}
\put( 0,25){Level $1$}
\put( 0,13){Level $2$}
\put( 0, 2){Level $3$}
\put(90,36){\footnotesize$I_{1} = [1,\,2,\,\dots,\,400]$}
\put(90,25){\footnotesize$I_{2} = [1,\,2,\,\dots,\,200]$, $I_{3} = [201,\,202,\,\dots,\,400]$}
\put(90,13){\footnotesize$I_{4} = [1,\,2,\,\dots,\,100]$, $I_{5} = [101,\,102,\,\dots,\,200]$, \dots}
\put(90, 2){\footnotesize$I_{8} = [1,\,2,\,\dots,\,50]$, $I_{9} = [51,\,52,\,\dots,\,100]$, \dots}
\put(52.5,36.5){$1$}
\put(35,25){$2$}
\put(70,25){$3$}
\put(26,13){$4$}
\put(44,13){$5$}
\put(61,13){$6$}
\put(78.5,13){$7$}
\put(22, 2){\small$8$}
\put(31, 2){\small$9$}
\put(38.5, 2){\small$10$}
\put(47, 2){\small$11$}
\put(56, 2){\small$12$}
\put(64.5, 2){\small$13$}
\put(73, 2){\small$14$}
\put(82, 2){\small$15$}
\end{picture}
\caption{Numbering of nodes in a fully populated binary tree with $L=3$ levels.
The root is the original index vector $I = I_{1} = [1,\,2,\,\dots,\,400]$.}
\label{fig:tree}
\end{figure}

The HBS representation
of an $M\times M$ matrix $\mtx{H}$ is
based on a partition of the index vector $I = [1,\,2,\,\dots,\,M]$
into a binary tree structure.
We let $I$ form the root of the tree, and give it the index $1$,
$I_{1} = I$. We next split the root into two roughly equi-sized
vectors $I_{2}$ and $I_{3}$ so that $I_{1} = I_{2} \cup I_{3}$.
The full tree is then formed by continuing to subdivide any interval
that holds more than some preset fixed number $m$ of indices.
We use the integers $\ell = 0,\,1,\,\dots,\,L$ to label the different
levels, with $0$ denoting the coarsest level.
A \textit{leaf} is a node corresponding to a vector that never got split.
For a non-leaf node $\tau$, its \textit{children} are the two boxes
$\sigma_{1}$ and $\sigma_{2}$ such that $I_{\tau} = I_{\sigma_{1}} \cup I_{\sigma_{2}}$,
and $\tau$ is then the \textit{parent} of $\sigma_{1}$ and $\sigma_{2}$.
Two boxes with the same parent are called \textit{siblings}. These
definitions are illustrated in Figure \ref{fig:tree}.

\subsection{Definition of the HBS property}
\label{sec:rigor}
We now define what it means for an $M\times M$ matrix $\mtx{H}$
to be \textit{hierarchically block separable} with respect to a given binary tree $\mathcal{T}$
that partitions the index vector $J = [1,\,2,\,\dots,\,M]$. For simplicity, we suppose that for every leaf node $\tau$ the index vector
$I_{\tau}$ holds precisely $m$ points, so that $M = m\,2^{L}$. Then $\mtx{H}$ is HBS with block
rank $k$ if the following two conditions hold:

\lsp

\noindent
\textit{(1) Assumption on ranks of off-diagonal blocks at the finest level:}
For any two distinct leaf nodes $\tau$ and $\tau'$, define the $n\times n$ matrix
\begin{equation}
\label{eq:def1}
\mtx{H}_{\tau,\tau'} = \mtx{H}(I_{\tau},I_{\tau'}).
\end{equation}
Then there must exist matrices $\mtx{U}_{\tau}$, $\mtx{V}_{\tau'}$, and $\tilde{\mtx{H}}_{\tau,\tau'}$ such that
\begin{equation}
\label{eq:def2}
\begin{array}{cccccccccc}
\mtx{H}_{\tau,\tau'} & = & \mtx{U}_{\tau} & \tilde{\mtx{H}}_{\tau,\tau'} & \mtx{V}_{\tau'}^{*}.  \\
m\times m &   & m\times k & k \times k & k\times m
\end{array}
\end{equation}

\lsp

\noindent
\textit{(2) Assumption on ranks of off-diagonal blocks on level $\ell = L-1,\,L-2,\,\dots,\,1$:}
The rank assumption at level $\ell$ is defined in terms of the blocks constructed
on the next finer level $\ell+1$: For any distinct nodes $\tau$ and
$\tau'$ on level $\ell$ with  children $\sigma_{1},\sigma_{2}$ and $\sigma_{1}',\sigma_{2}'$,
respectively, define
\begin{equation}
\label{eq:def3}
\mtx{H}_{\tau,\tau'} = \left[\begin{array}{cc}
\tilde{\mtx{H}}_{\sigma_{1},\sigma_{1}'} & \tilde{\mtx{H}}_{\sigma_{1},\sigma_{2}'} \\
\tilde{\mtx{H}}_{\sigma_{2},\sigma_{1}'} & \tilde{\mtx{H}}_{\sigma_{2},\sigma_{2}'}
\end{array}\right].
\end{equation}
Then there must exist matrices $\mtx{U}_{\tau}$, $\mtx{V}_{\tau'}$, and $\tilde{\mtx{H}}_{\tau,\tau'}$ such that
\begin{equation}
\label{eq:def4}
\begin{array}{cccccccccc}
\mtx{H}_{\tau,\tau'} & = & \mtx{U}_{\tau} & \tilde{\mtx{H}}_{\tau,\tau'} & \mtx{V}_{\tau'}^{*}.  \\
2k\times 2k &   & 2k\times k & k \times k & k\times 2k
\end{array}
\end{equation}

\lsp

 An HBS matrix is
now fully described if the basis matrices $\mtx{U}_{\tau}$ and $\mtx{V}_{\tau}$ are provided
for each node $\tau$, and in addition, we are for each leaf $\tau$ given the $m\times m$ matrix
\begin{equation}
\label{eq:def5}
\mtx{D}_{\tau} = \mtx{H}(I_{\tau},I_{\tau}),
\end{equation}
and for each parent node $\tau$ with children $\sigma_{1}$ and $\sigma_{2}$ we are given
the $2k\times 2k$ matrix
\begin{equation}
\label{eq:def6}
\mtx{B}_{\tau}       = \left[\begin{array}{cc}
0 & \tilde{\mtx{H}}_{\sigma_{1},\sigma_{2}} \\
\tilde{\mtx{H}}_{\sigma_{2},\sigma_{1}} & 0
\end{array}\right].
\end{equation}
Observe in particular that the matrices $\tilde{\mtx{H}}_{\sigma_{1},\sigma_{2}}$ are only required when
$\{\sigma_{1},\sigma_{2}\}$
forms a sibling pair. Figure \ref{fig:summary_of_factors} summarizes the required matrices.

\begin{figure}
\small
\begin{tabular}{|l|l|l|l|} \hline
                & Name:            & Size:       & Function: \\ \hline
For each leaf   & $\mtx{D}_{\tau}$ & $ m\times  m$ & The diagonal block $\mtx{H}(I_{\tau},I_{\tau})$. \\
node $\tau$:    & $\mtx{U}_{\tau}$ & $ m\times  k$ & Basis for the columns in the blocks in row $\tau$. \\
                & $\mtx{V}_{\tau}$ & $ m\times  k$ & Basis for the rows in the blocks in column $\tau$. \\ \hline
For each parent & $\mtx{B}_{\tau}$ & $2k\times 2k$ & Interactions between the children of $\tau$. \\
node $\tau$:    & $\mtx{U}_{\tau}$ & $2k\times  k$ & Basis for the columns in the (reduced) blocks in row $\tau$. \\
                & $\mtx{V}_{\tau}$ & $2k\times  k$ & Basis for the rows in the (reduced) blocks in column $\tau$. \\ \hline
\end{tabular}
\caption{An HBS matrix $\mtx{H}$ associated with a tree $\mathcal{T}$ is fully specified if the factors
listed above are provided.}
\label{fig:summary_of_factors}
\end{figure}

\subsection{Telescoping factorization}
Given the matrices defined in the previous section,  we define the following block diagonal factors:
\begin{align}
\label{eq:def_uD}
\underline{\mtx{D}}^{(\ell)} &= \mbox{diag}(\mtx{D}_{\tau}\,\colon\, \tau\mbox{ is a box on level }\ell),\qquad \ell = 0,\,1,\,\dots,\,L,\\
\underline{\mtx{U}}^{(\ell)} &= \mbox{diag}(\mtx{U}_{\tau}\,\colon\, \tau\mbox{ is a box on level }\ell),\qquad \ell = 1,\,2,\,\dots,\,L,\\
\underline{\mtx{V}}^{(\ell)} &= \mbox{diag}(\mtx{V}_{\tau}\,\colon\, \tau\mbox{ is a box on level }\ell),\qquad \ell = 1,\,2,\,\dots,\,L,\\
\underline{\mtx{B}}^{(\ell)} &= \mbox{diag}(\mtx{B}_{\tau}\,\colon\, \tau\mbox{ is a box on level }\ell),\qquad \ell = 0,\,1,\,\dots,\,L-1,.
\end{align}
 Furthermore, we let $\tilde{\mtx{H}}^{(\ell)}$ denote the block matrix whose diagonal blocks are zero,
and whose off-diagonal blocks are the blocks $\tilde{\mtx{H}}_{\tau,\tau'}$ for all distinct $\tau,\tau'$
on level $\ell$.
With these definitions,
\begin{equation}
\label{eq:tele1}
\begin{array}{cccccccccccc}
\mtx{H} & = & \underline{\mtx{U}}^{(L)} & \tilde{\mtx{H}}^{(L)} & (\underline{\mtx{V}}^{(L)})^{*} & + & \underline{\mtx{D}}^{(L)};\\
m\,2^{L} \times n\,2^{L} &&
m\,2^{L} \times k\,2^{L} &
k\,2^{L} \times k\,2^{L} &
k\,2^{L} \times m\,2^{L} &&
m\,2^{L} \times m\,2^{L}
\end{array}
\end{equation}
for $\ell = L-1,\,L-2,\,\dots,\,1$ we have
\begin{equation}
\label{eq:tele2}
\begin{array}{cccccccccccc}
\tilde{\mtx{H}}^{(\ell+1)} &=& \underline{\mtx{U}}^{(\ell)} & \tilde{\mtx{H}}^{(\ell)} & (\underline{\mtx{V}}^{(\ell)})^{*} & +
& \underline{\mtx{B}}^{(\ell)};\\
k\,2^{\ell+1} \times k\,2^{\ell+1} &&
k\,2^{\ell+1} \times k\,2^{\ell} &
k\,2^{\ell} \times k\,2^{\ell} &
k\,2^{\ell} \times k\,2^{\ell+1} &&
k\,2^{\ell+1} \times k\,2^{\ell+1}
\end{array}
\end{equation}
and finally
\begin{equation}
\label{eq:tele3}
\tilde{\mtx{H}}^{(1)} = \underline{\mtx{B}}^{(0)}.
\end{equation}

\section{Fast arithmetic operations on HBS matrices}
\label{sec:fastHBS}

Arithmetic operations involving dense HBS matrices of size $N\times N$ can often
be executed in $O(N)$ operations. This fast matrix algebra is vital for achieving
linear complexity in our direct solver. This section provides a brief introduction
to the HBS matrix algebra. We describe the operations we need (inversion, addition,
and low-rank update) in some detail for the single level ``block separable'' format.
The generalization to the multi-level ``hierarchically block separable'' format is
briefly described for the case of matrix inversion. A full description of all
algorithms required is given in \cite{Adiss}, which is related to the earlier work
\cite{2010_jianlin_fast_hss}.

Before we start, we recall that a block separable matrix $\mtx{H}$ consisting of
$p\times p$ blocks, each of size $m\times m$, and with ``HBS-rank'' $k<m$, admits
the factorization
\begin{equation}
\label{eq:yy2b}
\begin{array}{cccccccccc}
 \mtx{H} &  =& \mtx{U}&\tilde{\mtx{H}}& \mtx{V}^{*} & +&  \mtx{D}.\\
mp\times mp &   & mp\times kp & kp \times kp & kp\times mp && mp \times mp\\
\end{array}
\end{equation}

\subsection{Inversion of a block separable matrix}
The decomposition (\ref{eq:yy2b}) represents $\mtx{H}$ as a sum of
one term $\mtx{U}\tilde{\mtx{H}}\mtx{V}^{*}$ that is ``low rank,''
and one term $\mtx{D}$ that is easily invertible (since it is block diagonal).
By modifying the classical Woodbury formula for inversion of a matrix
perturbed by the addition of a low-rank term, it can be shown that (see Lemma 3.1 of
\cite{m2011_1D_survey})
\begin{equation}
\label{eq:woodbury}
\mtx{H}^{-1} = \mtx{E}\,(\tilde{\mtx{H}} + \hat{\mtx{D}})^{-1}\,\mtx{F}^{*} + \mtx{G},
\end{equation}
where
\begin{align}
\label{eq:def_muhD}
\hat{\mtx{D}} =&\ \bigl(\mtx{V}^{*}\,\mtx{D}^{-1}\,\mtx{U}\bigr)^{-1},\\
\label{eq:def_muE}
\mtx{E}  =&\ \mtx{D}^{-1}\,\mtx{U}\,\hat{\mtx{D}},\\
\label{eq:def_muF}
\mtx{F}  =&\ (\hat{\mtx{D}}\,\mtx{V}^{*}\,\mtx{D}^{-1})^{*},\\
\label{eq:def_muG}
\mtx{G}  =&\ \mtx{D}^{-1} - \mtx{D}^{-1}\,\mtx{U}\,\hat{\mtx{D}}\,\mtx{V}^{*}\,\mtx{D}^{-1},
\end{align}
assuming the inverses in formulas (\ref{eq:woodbury}) --- (\ref{eq:def_muG}) all exist.
Now observe that the matrices $\hat{\mtx{D}}$, $\mtx{E}$, $\mtx{F}$, and $\mtx{G}$
can all easily be computed since the formulas defining them involve only block-diagonal matrices.
In consequence, (\ref{eq:woodbury}) reduces the task of inverting the big (size $mp\times mp$)
matrix $\mtx{H}$ to the task of inverting the small (size $kp\times kp$) matrix $\tilde{\mtx{H}} + \hat{\mtx{D}}$.

When $\mtx{H}$ is not only ``block separable'', but ``hierarchically block separable'', the
process can be repeated recursively by exploiting that $\tilde{\mtx{H}} + \hat{\mtx{D}}$ is
itself amenable to accelerated inversion, etc. The resulting process is somewhat tricky to
analyze, but leads to very clean codes. To illustrate, we include Algorithm 1
which shows the multi-level $O(N)$ inversion algorithm for an HBS matrix $\mtx{H}$.
The algorithm takes as input the factors $\{\mtx{U}_{\tau},\,\mtx{V}_{\tau},\,\mtx{D}_{\tau},\,\mtx{B}_{\tau}\}_{\tau}$
representing $\mtx{H}$ (cf.~Figure \ref{fig:summary_of_factors}), and outputs an analogous
set of factors $\{\mtx{E}_{\tau},\,\mtx{F}_{\tau},\,\mtx{G}_{\tau}\}_{\tau}$ representing
$\mtx{H}^{-1}$. With these factors, the matrix-vector multiplication $\vct{y} = \mtx{H}^{-1}\vct{x}$
can be executed via the procedure described in Algorithm 2.

\begin{figure}[ht]
\begin{center}
\fbox{
\begin{minipage}{.9\textwidth}\small

\begin{center}
\textsc{Algorithm 1} (inversion of an HBS matrix)
\end{center}

Given factors $\{\mtx{U}_{\tau},\,\mtx{V}_{\tau},\,\mtx{D}_{\tau},\,\mtx{B}_{\tau}\}_{\tau}$
representing an HBS matrix $\mtx{H}$, this algorithm constructs factors
$\{\mtx{E}_{\tau},\,\mtx{F}_{\tau},\,\mtx{G}_{\tau}\}_{\tau}$ representing
$\mtx{H}^{-1}$.

\begin{tabbing}
\hspace{5mm} \= \hspace{5mm} \= \hspace{5mm} \= \kill
\textbf{loop} over all levels, finer to coarser, $\ell = L,\,L-1,\,\dots,1$\\
\> \textbf{loop} over all boxes $\tau$ on level $\ell$,\\
\> \> \textbf{if} $\tau$ is a leaf node\\
\> \> \> $\tilde{\mtx{D}}_{\tau} = \mtx{D}_{\tau}$\\
\> \> \textbf{else}\\
\> \> \> Let $\sigma_{1}$ and $\sigma_{2}$ denote the children of $\tau$.\\
\> \> \> $\tilde{\mtx{D}}_{\tau} = \mtwo{\hat{\mtx{D}}_{\sigma_{1}}}{\mtx{B}_{\sigma_{1},\sigma_{2}}}{\mtx{B}_{\sigma_{2},\sigma_{1}}}{\hat{\mtx{D}}_{\sigma_{2}}}$\\
\> \> \textbf{end if}\\
\> \> $\hat{\mtx{D}}_{\tau} = \bigl(\mtx{V}_{\tau}^{*}\,\tilde{\mtx{D}}_{\tau}^{-1}\,\mtx{U}_{\tau}\bigr)^{-1}$.\\
\> \> $\mtx{E}_{\tau} = \tilde{\mtx{D}}_{\tau}^{-1}\,\mtx{U}_{\tau}\,\hat{\mtx{D}}_{\tau}$.\\
\> \> $\mtx{F}_{\tau}^{*} = \hat{\mtx{D}}_{\tau}\,\mtx{V}_{\tau}^{*}\,\tilde{\mtx{D}}_{\tau}^{-1}$.\\
\> \> $\mtx{G}_{\tau} = \tilde{\mtx{D}}_{\tau}^{-1} - \tilde{\mtx{D}}_{\tau}^{-1}\,\mtx{U}_{\tau}\,\hat{\mtx{D}}_{\tau}\,\mtx{V}_{\tau}^{*}\,\tilde{\mtx{D}}_{\tau}^{-1}$.\\
\> \textbf{end loop}\\
\textbf{end loop}\\
$\mtx{G}_{1} = \mtwo{\hat{\mtx{D}}_{2}}{\mtx{B}_{2,3}}{\mtx{B}_{3,2}}{\hat{\mtx{D}}_{3}}^{-1}$.
\end{tabbing}
\end{minipage}}
\end{center}
\end{figure}

\begin{figure}
\begin{center}
\fbox{
\begin{minipage}{.9\textwidth}\small

\begin{center}
\textsc{Algorithm 2} (application of inverse)
\end{center}

\textit{Given $\vct{x}$, compute $\vct{y} = \mtx{H}^{-1}\,\vct{x}$ using the
compressed representation of $\mtx{H}^{-1}$ resulting from Algorithm 1.}

\begin{tabbing}
\hspace{5mm} \= \hspace{5mm} \= \hspace{5mm} \= \kill
\textbf{loop} over all leaf boxes $\tau$\\
\> $\hat{\vct{x}}_{\tau} = \mtx{F}_{\tau}^{*}\,\vct{x}(I_{\tau})$.\\
\textbf{end loop}\\[1mm]
\textbf{loop} over all levels, finer to coarser, $\ell = L,\,L-1,\,\dots,1$\\
\> \textbf{loop} over all parent boxes $\tau$ on level $\ell$,\\
\> \> Let $\sigma_{1}$ and $\sigma_{2}$ denote the children of $\tau$.\\
\> \> $\hat{\vct{x}}_{\tau} = \mtx{F}_{\tau}^{*}\,\vtwo{\hat{\vct{x}}_{\sigma_{1}}}{\hat{\vct{x}}_{\sigma_{2}}}$.\\
\> \textbf{end loop}\\
\textbf{end loop}\\[1mm]
$\vtwo{\hat{\vct{y}}_{2}}{\hat{\vct{y}}_{3}} = \mtx{\mtx{G}}_{1}\,\vtwo{\hat{\vct{x}}_{2}}{\hat{\vct{x}}_{3}}$.\\[1mm]
\textbf{loop} over all levels, coarser to finer, $\ell = 1,\,2,\,\dots,\,L-1$\\
\> \textbf{loop} over all parent boxes $\tau$ on level $\ell$\\
\> \> Let $\sigma_{1}$ and $\sigma_{2}$ denote the children of $\tau$.\\
\> \> $\vtwo{\hat{\vct{y}}_{\sigma_{1}}}{\hat{\vct{y}}_{\sigma_{2}}} =
       \mtx{E}_{\tau}\,\hat{\vct{x}}_{\tau} +
       \mtx{G}_{\tau}\,\vtwo{\hat{\vct{x}}_{\sigma_{1}}}{\hat{\vct{x}}_{\sigma_{2}}}$.\\
\> \textbf{end loop}\\
\textbf{end loop}\\[1mm]
\textbf{loop} over all leaf boxes $\tau$\\
\> $\vct{y}(I_{\tau}) = \mtx{E}_{\tau}\,\hat{\vct{q}}_{\tau} + \mtx{G}_{\tau}\,\vct{x}(I_{\tau})$.\\
\textbf{end loop}
\end{tabbing}
\end{minipage}}
\end{center}
\end{figure}


\subsection{Addition of two block separable matrices}
Let $\mtx{H}^A$ and $\mtx{H}^B$ be block separable matrices with factorizations
$$
\mtx{H}^A   = \mtx{U}^A \tilde{\mtx{H}}^A \mtx{V}^{A*}  +  \mtx{D}^A,
\qquad\mbox{and}\qquad
\mtx{H}^B   = \mtx{U}^B \tilde{\mtx{H}}^B \mtx{V}^{B*}  +  \mtx{D}^B.
$$
Then $\mtx{H} = \mtx{H}^{A}+\mtx{H}^{B}$ can be written in block separable form via
\begin{equation}
\label{eq:postadd}
\mtx{H} = \mtx{H}^{A}+\mtx{H}^{B} = \left[\mtx{U}^A\, \mtx{U}^B\right]
\mtwo{\tilde{\mtx{H}}^A}{0}{0}{\tilde{\mtx{H}}^B} \left[\mtx{V}^A\, \mtx{V}^B\right]^*
+\left(\mtx{D}^A+\mtx{D}^B\right).
\end{equation}
To restore (\ref{eq:postadd}) to block separable form, permute the rows and columns
of $\left[\mtx{U}^A\, \mtx{U}^B\right]$ and $\left[\mtx{V}^A\, \mtx{V}^B\right]$ to
attain block diagonal form, then re-orthogonalize the diagonal blocks. This process
in principle results in a matrix $\mtx{H}$ whose HBS-rank is the sum of the HBS-ranks
of $\mtx{H}^{A}$ and $\mtx{H}^{B}$. In practice, this rank increase can be combatted
by numerically recompressing the basis matrices, and updating the middle factor as
needed. For details, as well as the extension to a multi-level scheme, see
\cite{2010_jianlin_fast_hss,Adiss}.

\subsection{Addition of a block separable matrix with a low rank matrix}
Let $\mtx{H}^{B}= \mtx{Q}\mtx{R}$ be a $k$-rank matrix where $\mtx{Q}$ and $\mtx{R}^*$ are of size $mp\times k$.  We
would like to add $\mtx{H}^B$ to the block separable matrix $\mtx{H}^A$.
Since we already know how to add two block separable matrices, we choose
to rewrite $\mtx{H}^B$ in block separable form.  Without
loss of generality, assume $\mtx{Q}$ is orthogonal.  Partition $\mtx{Q}$ into $p$ blocks of size $m\times k$.
The blocks make up the matrix $\mtx{U}^B$.  Likewise partition $\mtx{R}$ into $p$ blocks of size $k \times m$.
The block matrix $\mtx{D}^B$ has entries $\mtx{D}_{\tau} = \mtx{Q}_{\tau}\mtx{R}_{\tau}$ for $\tau = 1,\ldots,p$.
  To construct the matrices $\mtx{V}^B$, for each $\tau = 1,\ldots,p$, the matrix $\mtx{R}_{\tau}$ is
factorized into $\tilde{\mtx{R}}_{\tau} \mtx{V}_{\tau}*$ where the matrix $\mtx{V}_{\tau}$ is
orthogonal. The matrices $\tilde{\mtx{R}}_{\tau}$ make up the entries of $\tilde{\mtx{H}}^B$.

\section{Accelerating the nested dissection algorithm}
\label{sec:accel}
In this section, we apply the structured matrix techniques introduced in
Sections \ref{sec:HBS} and \ref{sec:fastHBS} to reduce the complexity of the solver
of Section \ref{sec:merge} from $O(N^{1.5})$ to $O(N)$. The key task that
we need to accelerate is the construction of the Schur complement for a
parent box from the Schur complements of its two children. The formula that
needs to be evaluated is, cf.~(\ref{eq:S_fromtwo_nd}),
\begin{equation}
\label{eq:napoli}
\mtx{S} =
\left[\begin{array}{cc}
\mtx{S}_{11} & \mtx{A}_{12} \\
\mtx{A}_{21} & \mtx{S}_{22}
\end{array}\right] -
\left[\begin{array}{cc}
\mtx{S}_{13} & \mtx{A}_{14}    \\
\mtx{A}_{23}     & \mtx{S}_{24}
\end{array}\right]\,
\left[\begin{array}{cc}
\mtx{S}_{33} & \mtx{A}_{34} \\
\mtx{A}_{43} & \mtx{S}_{44}
\end{array}\right]^{-1}\,
\left[\begin{array}{cc}
\mtx{S}_{31} & \mtx{A}_{32}      \\
\mtx{A}_{41}      & \mtx{S}_{42}
\end{array}\right].
\end{equation}

The acceleration can be broken into three steps which utilize important properties of each submatrix.

\begin{itemize}
 \item[\textbf{Step 1}]  The inverse in equation (\ref{eq:napoli}) never needs to be constructed.  Instead the solution of
\begin{equation}
\label{eq:step1}
\left[\begin{array}{cc}
\mtx{S}_{33} & \mtx{A}_{34} \\
\mtx{A}_{43} & \mtx{S}_{44}
\end{array}\right]\,\vtwo{\mtx{X}_3}{\mtx{X}_4} = \vtwo{\mtx{Z}_3}{\mtx{Z}_4},
\end{equation}
can be found by rapidly via a block solve.
Then $\vtwo{\mtx{X}_3}{\mtx{X}_4}$ is given by
$$\mtx{X}_4 = \left(\mtx{S}_{44}-\mtx{A}_{43}\mtx{S}_{33}^{-1}\mtx{A}_{34}\right)^{-1}\left(\mtx{Z}_4-
\mtx{A}_{43}\mtx{S}_{33}^{-1}\mtx{Z}_3\right)$$
and
$$\mtx{X}_3 = \mtx{S}_{33}^{-1}\mtx{Z}_3-\mtx{S}_{33}^{-1}\mtx{A}_{34}\mtx{X}_4.$$

Since $\mtx{S}_{33}$ is HBS, an approximation of its inverse can be constructed and applied rapidly.
The matrix $\mtx{A}_{43}\mtx{S}_{33}^{-1}\mtx{A}_{34}$
is also HBS, since $\mtx{A}_{43}$ and $\mtx{A}_{34}$ are anti-diagonal matrices (ie. all the entries are zero
except those on the diagonal going from the lower left corner to the upper right corner).
Hence $\mtx{S}_{44}-\mtx{A}_{43}\mtx{S}_{33}^{-1}\mtx{A}_{34}$ can be
added quickly.  The resulting matrix is HBS and can be
inverted with linear scaling computational cost.

Let $\left[\begin{array}{cc}
\mtx{X}_{31} & \mtx{X}_{32}      \\
\mtx{X}_{41}      & \mtx{X}_{42}
\end{array}\right]$ denote the result of applying the block inverse to $\left[\begin{array}{cc}
\mtx{S}_{31} & \mtx{A}_{32}      \\
\mtx{A}_{41}      & \mtx{S}_{42}
\end{array}\right]$.

\item[\textbf{Step 2}]  The matrices $\mtx{S}_{13}$, and $\mtx{S}_{24}$ are low rank, thus we can rewrite the matrices
in their low rank factored form as $\mtx{L}_{13}\mtx{R}_{13}$ and $\mtx{L}_{24}\mtx{R}_{24}$.
Using this notation, the second term in (\ref{eq:napoli}) can be expressed in a low rank factored form
\begin{equation}
\left[\begin{array}{cc}
\mtx{S}_{13} & \mtx{A}_{14}    \\
\mtx{A}_{23}     & \mtx{S}_{24}
\end{array}\right]\, \left[\begin{array}{cc}
\mtx{X}_{31} & \mtx{X}_{32}      \\
\mtx{X}_{41}      & \mtx{X}_{42}
\end{array}\right]  = \left[\begin{array}{cc}
\mtx{L}_{13}\mtx{R}_{13}\mtx{X}_{31}+\mtx{A}_{14}\mtx{X}_{41} & \mtx{L}_{13}\mtx{R}_{13}\mtx{X}_{32}  +\mtx{A}_{14}\mtx{X}_{42}    \\
\mtx{A}_{23}\mtx{X}_{31}+\mtx{L}_{24}\mtx{R}_{24}\mtx{X}_{41} & \mtx{A}_{23}\mtx{X}_{32}+\mtx{L}_{24}\mtx{R}_{24}\mtx{X}_{42}
\end{array}\right].\label{eq:matmult}\end{equation}
Since all the matrices $\mtx{A}_{jk}$ are very sparse, the four blocks on the right hand side of
(\ref{eq:matmult}) are of low rank.  (Recall in the case of the
five point stencil, the matrices $\mtx{A}_{14}$ and $\mtx{A}_{23}$ are zero).

  Let $\left[\begin{array}{cc}
\mtx{L}_{11}\mtx{R}_{11} & \mtx{L}_{12}\mtx{R}_{12}    \\
\mtx{L}_{21}\mtx{R}_{21} & \mtx{L}_{22}\mtx{R}_{22}
\end{array}\right]$ denote the low rank factorization of the blocks in (\ref{eq:matmult}).

\item[\textbf{Step 3}] Now we add the two terms that comprise the Schur complement
$$\left[\begin{array}{cc}
\mtx{S}_{11} & \mtx{A}_{12} \\
\mtx{A}_{21} & \mtx{S}_{22}
\end{array}\right] -\left[\begin{array}{cc}
\mtx{L}_{11}\mtx{R}_{11} & \mtx{L}_{12}\mtx{R}_{12}    \\
\mtx{L}_{21}\mtx{R}_{21} & \mtx{L}_{22}\mtx{R}_{22}
\end{array}\right].$$
The diagonal block entries are HBS + low rank which is computed via Algorithms \ref{alg:QR_HSS} and \ref{alg:HSS_add}.  The off-diagonal
blocks are low rank with a very sparse update which is also low rank.  The result is one HBS matrix.

\end{itemize}

\begin{remark}
As a practical matter, structured matrix algebra should not be introduced
until the Schur complements get fairly large (roughly of size $1000 \times 1000$ or so). 
This means that at the lower levels, formula (\ref{eq:napoli}) is evaluated using
dense matrix algebra for all matrices $\mtx{S}_{i,j}$.
\end{remark}

\section{Numerical experiments}
\label{sec:num}
In this section, we illustrate the capabilities of the proposed method for constructing
solution operators for problems of the form
\begin{equation}
\label{eq:lap_psbc}
\left\{\begin{aligned}
- \Delta u(\vct{x}) + b(\vct{x})u_x (\vct{x}) + c(\vct{x})u_y (\vct{x})+d(\vct{x})u(\vct{x}) & = f(\vct{x}),
\quad  &&\vct{x} \in \Omega =[0,1]^2,\\
         u(\vct{x}) &= g(\vct{x}),       \quad    &&\vct{x} \in \Gamma,
\end{aligned}\right.
\end{equation}
where $b$, $c$, $d$, and $f$ are functions defined on $\Omega$, and the boundary data $g$ is defined on $\Gamma$.
Section \ref{sec:explore} investigates the asymptotic complexity of the direct solver for several
different differential operators (Laplace, Helmholtz, convection-diffusion, etc.) for the
case where the body load is zero ($f=0$). It also reports on the accuracy for each case.
Section \ref{sec:scale} reports the execution times of the build and the solve stages of the direct solver.
Section \ref{sec:bodyload} reports on the performance of the method for a problem with localized body loads.

For all problems, the domain is discretized with a uniform grid of $n \times n$ points so that the
grid spacing is $h = 1/(n-1)$. We let $N = n^{2}$ denote the total number of nodes.
Equation (\ref{eq:lap_psbc}) is discretized with the finite difference scheme
corresponding to the five point stencil.  For example, when a node $k$ is in the interior of $\Omega$, the
discretization of the differential operator in (\ref{eq:lap_psbc}) is
\begin{multline*}
\frac{1}{h^2}\bigl[4u(k)-u(k_n)-u(k_s)-u(k_w)-u(k_e)\bigr] +\\
\frac{1}{h}b(k)\bigl[u(k_e)-u(k_w)\bigr]+
\frac{1}{h}c(k)\bigl[u(k_n)-u(k_s)\bigr]+d(k)u(k),
\end{multline*}
where $k_e$, $k_n$, $k_w$, $k_s$ denote the grid points to the ``east'', ``north'', ``west'', and
``south'' of $k$, respectively. The HBS matrix algebra was run with a local tolerance
of $\varepsilon = 10^{-7}$.

All experiments are executed on a Lenovo laptop computer with a 2.4GHz Intel i5 processor and 8GB of RAM. 
The method was implemented in Matlab, which we judged adequate since the main purpose of the
experiments is to substantiate our claims in regards to asymptotic complexity.
It should be noted, however, that even this non-optimized code runs quite fast, see, e.g., Table \ref{tab:nested_times}.


Recall that the approximate solution (or Dirichlet-to-Neumann) operator $\mtx{G}$ is the inverse
of the Schur complement $\mtx{S}$ for the full domain, see Section \ref{sec:schur}.

\subsection{Range of problems with optimal scaling }
\label{sec:explore}
The proposed method for constructing Dirichlet-to-Neumann operators has been applied
to several problems to investigate its asymptotic complexity. The problems are:
\begin{itemize}
\item \textit{Laplace:}  Let $b(\vct{x}) = c(\vct{x}) = d(\vct{x})=0$.
\item \textit{Diffusion-Convection I:}
Let $c(\vct{x})=d(\vct{x})=0$ and the convection in the $x$ direction be constant: $b(\vct{x}) = 100$.
\item \textit{Diffusion-Convection II:} Same as \textit{Diffusion-Convection I}, but with $b(\vct{x}) = 1000$.
\item \textit{Diffusion-convection III:} Introduce a divergence free convection field by setting
$b(\vct{x})= 125 \cos(4 \pi y )$ and $c(\vct{x}) = 125 \sin(4 \pi x )$, and $d(\vct{x})=0$.
\item \textit{Diffusion-convection IV:} Introduce a convection field with sources and sinks by setting
 Let $b(\vct{x})= 125 \cos(4 \pi x )$, $c(\vct{x}) = 125 \sin(4 \pi y )$, and $d(\vct{x})=0$.
\item \textit{Helmholtz I:} Consider the Helmholtz equation corresponding to a domain
that is roughly $1.5\times 1.5$ wavelengths large: $b(\vct{x}) = c(\vct{x})=0$ and $d(\vct{x}) = -100$.
\item \textit{Helmholtz II:} Consider the Helmholtz equation corresponding to a domain
that is roughly $10\times 10$ wavelengths large: $b(\vct{x}) = c(\vct{x})=0$ and $d(\vct{x}) = -4005$.
\item \textit{Helmholtz III:} Consider the Helmholtz equation near a resonance:
$b(\vct{x}) = c(\vct{x})=0$ and $d(\vct{x}) = -\lambda_{10} +10^{-5}$, where
$\lambda_{10}$ is the tenth eigenvalue of the discrete Laplace operator (note that these are known analytically).
\item \textit{Helmholtz IV:} Consider a sequence of Helmholtz problems where the wave-number
is increased to keep a constant 40 points per wave-length:
$b(\vct{x}) = c(\vct{x})=0$ and $d(\vct{x}) = -\left(\frac{2\pi n}{40}\right)^2$.
\item \textit{Random Laplacian I:} Let the matrix $\mtx{A}$ reflect an elliptic
equilibrium problem on a network instead of a continuum PDE.  In this case, the
network is the square grid where each link is assigned a random \textit{conductivity}
between varying between $1$ and $2$. The potential at any single node is the weighted average
of the potentials of its four neighbors, where the weights are the conductivities.
\item \textit{Random Laplacian II:} Same as \textit{Random Laplacian I}, but now the
conductivities vary between $1$ and $1000$.
\end{itemize}

Table \ref{tab:mem} reports the amount of memory $M(n)$ in MB required to store the
Dirichlet-to-Neumann operator $\mtx{G}$ for each problem; it also reports the fraction $M(n)/n$.
Our claim in regards to compressibility amounts to a prediction that $M(n)/n$ will remain
stable as $n$ grows for all problems, except for Helmholtz IV.
Table \ref{tab:mem} demonstrates that this scaling holds true in the range $n = [256,\,512,\,1024,\,2048]$.

Table \ref{tab:error} reports two errors measured on a grid of size $1024\times 1024$:\\
\begin{tabular}{r l}
$e_1$ &- the relative $l^{2}$-error in the vector $\mtx{G}\,\vct{r}$ where $\vct{r}$ is a unit vector of random direction\\
$e_2$ &- the relative $l^{2}$-error in the vector $\mtx{G}\,\vct{r}$ where $\vct{r}$ is smooth.
\end{tabular}\\
The exact value of $\mtx{G}\,\vct{r}$ was found by using GMRES to solve the full original
linear system $\mtx{A}\vct{x} = \hat{\vct{r}}$, where $\hat{\vct{r}}$ is a vector of length
$n^{2}$ such that $\hat{\vct{r}}|_{\Gamma} = \vct{r}$ and $\hat{\vct{r}}|_{\Omega/\Gamma}  = 0$.
A slight loss in accuracy is observed for Helmholtz I, IV and Random II problems.  There is a
substantial loss in accuracy for the Helmholtz III problem.  This is to be expected since the matrix
$\mtx{A}$ is close to being numerically singular.

\begin{table}[ht]
\centering
\begin{tabular}{|c|c|c|c|c| }
\hline
  \textit{Problem}       & $n = 256$ & $n = 512$& $n = 1024$& $n = 2048$ \\ \hline
  \textit{Laplace}       &  0.83  (3.2e-3) &  1.62  (3.2e-3) &  3.18 (3.1e-3) &  6.27 (3.1e-3) \\ \hline
  \textit{DiffConv I}    &  0.91 (3.5e-3)  &  1.75 (3.4e-3)  &  3.32 (3.2e-3)  &  6.52 (3.2e-3) \\ \hline
  \textit{DiffConv II}   &  1.10 (4.3e-3)  & 1.84 (3.6e-3)  &  3.62 (3.5e-3)  &  6.87 (3.4e-3) \\ \hline
  \textit{DiffConv III}  &  0.86 (3.4e-3)  &  1.70 (3.3e-3)  &  3.32 (3.2e-3)  &  6.55 (3.3e-3) \\ \hline
  \textit{DiffConv IV}   &  0.97 (3.8e-3)  &  1.83 (3.6e-3)  &  3.43 (3.3e-3)  &  6.59 (3.2e-3) \\ \hline
  \textit{Helmholtz I}   &  0.86 (3.4e-3)  &  1.67 (3.3e-3)  &  3.25 (3.2e-3)  &  6.34 (3.1e-3) \\ \hline
  \textit{Helmholtz II}  &  1.04 (4.1e-3)  &  1.91 (3.7e-3)  &  3.56 (3.5e-3)  &  6.78 (3.3e-3) \\ \hline
  \textit{Helmholtz III} &  0.86 (3.4e-3)  &  1.67 (3.3e-3)  &  3.29 (3.2e-3)  & 6.42 (3.1e-3) \\ \hline
  \textit{Helmholtz IV}  &  0.89 (3.5e-3)  &  1.74 (3.4e-3)  &  3.59 (3.5e-3)  &  7.89 (3.9e-3) \\ \hline
  \textit{Random I}      &  0.83 (3.2e-3)  &  1.64 (3.2e-3)  &  3.22 (3.1e-3)  & 6.34 (3.1e-3) \\ \hline
  \textit{Random II}     &  0.82 (3.2e-3)  &  1.64 (3.2e-3)  &  3.23 (3.2e-3)  & 6.36 (3.1e-3) \\ \hline
\end{tabular}
\caption{\label{tab:mem} Memory $M(n)$ in MB required to store the solution operator
for the problems listed in Section \ref{sec:explore}. The quantity $\frac{M(n)}{n}$ is reported in parenthesis. }
\end{table}

\begin{table}[ht]
\centering
\begin{tabular}{|c|c|c| }
\hline
  \textit{Problem}       & $e_1$ & $e_2$ \\ \hline
  \textit{Laplace}       &  6.3e-7  & 3.6e-7 \\ \hline
  \textit{DiffConv I}    &  1.5e-6  &  1.3e-6\\ \hline
  \textit{DiffConv II}   &  8.7e-6  & 8.2e-6 \\ \hline
  \textit{DiffConv III}  & 5.6e-7 & 3.4e-7 \\ \hline
  \textit{DiffConv IV}   &  4.1e-8  & 4.1e-8 \\ \hline
  \textit{Helmholtz I}   &  1.4e-4  &  4.8e-7 \\ \hline
  \textit{Helmholtz II}  & 1.1e-6  &  5.1e-6 \\ \hline
  \textit{Helmholtz III} &  1.2e-5  &  5.7e-4\\ \hline
  \textit{Helmholtz IV}  &  8.2e-4  &  1.2e-3 \\ \hline
  \textit{Random I}      &  1.8e-7  &  1.2e-7 \\ \hline
  \textit{Random II}     &  1.4e-5  &  8.1e-6 \\ \hline
\end{tabular}
\caption{\label{tab:error} Errors $e_1$ and $e_2$ for the solution operator
for the problems listed in Section \ref{sec:explore}.}
\end{table}

\subsection{Performance}
\label{sec:scale}

In this section we report the computational times required for the
experiments that were described in Section \ref{sec:explore},
for grids of size $512 \times 512$ to $4096 \times 4096$.
Since the times were very similar across most experiments, we report only those for
\textit{Laplace}, which represents the ``typical'' times observed, and for \textit{Helmholtz IV},
which was the most challenging and slowest of all the experiments conducted.
Table \ref{tab:nested_times} reports:\\
\begin{tabular}{r l}
$T_{\rm build}$ &- the time in seconds for constructing the Dirichlet-to-Neumann operator,\\
$T_{\rm solve}$ &- the time in seconds for applying the Dirichlet-to-Neumann operator to a vector.\\
\end{tabular}\\
Let us draw the reader's attention to some interesting results in Table \ref{tab:nested_times}:
\begin{itemize}
\item A principal claim we make in terms of performance of the proposed method is that
after an initial pre-computation in which the solution operator is built, the time $T_{\rm solve}$
required to process a new vector of Dirichlet data is small. Table \ref{tab:nested_times} clearly
bears out our claim that $T_{\rm solve}$ scales linearly with \textit{the number of points on the
boundary}, in other words $T_{\rm solve} \sim N^{0.5}$. Moreover, the scaling constant turns out
to be small, for a grid of size $4096 \times 4096$, the solve time is only 0.1 seconds.
\item The other key claim made is that the time to build the solution operator in the first place
scales linearly with the number of points in the grid, in other words $T_{\rm build} \sim N$.
Table \ref{tab:nested_times} shows that for grids holding between 250k and 16M nodes, the build
time in fact scales sub-linearly. Eventually, linear complexity must of course kick in, but it is
interesting that it has not yet done so even for a grid holding over 16M nodes.
\item For the example \textit{Helmholtz IV} we did not predict linear complexity. This problem
models wave-propagation in such a way that as $n$ grows, the number of wave-lengths along a side
of the domain grows proportionally. This will eventually destroy the rank-structure in the Schur
complements that we rely on to reduce the $O(N^{1.5})$ scaling of classical nested dissection down
to $O(N)$. Now what is interesting is that while the \textit{predicted} complexity is
$T_{\rm build} \sim N^{1.5}$, the \textit{observed} complexity is only $T_{\rm build} \sim N$.
We expect that the predicted asymptotic scaling will eventually assert itself, but it is to
us remarkable that it has not yet done so given that the largest domain with 16M nodes represents
a physical problem of size $100 \times 100$ wavelengths.
\end{itemize}

\begin{table}[ht]
\centering
 \begin{tabular}{|c|c|c|c|c|}
\hline
   $N$ & \multicolumn{2}{|c|}{\textit{Laplace}}& \multicolumn{2}{|c|}{\textit{Helmholtz IV}}\\ \hline
&$T_{\rm build}$& $T_{\rm solve}$&$T_{\rm build}$& $T_{\rm solve}$\\
   &  (sec) & (sec)    &  (sec) & (sec)\\ \hline
  $512^2$ & $13.44$ & $0.013$ & $50.78$ & $0.013$\\ \hline
  $1024^2$ & $45.25$ & $0.027$ &$193.58$ &$0.027$\\ \hline
 $2048^2$ & $135.01$ & $0.058$ &$765.35$ & $0.056$ \\ \hline
 $4096^2$ & $450.73$ & $0.107$  & $3167.56$& $0.115$\\ \hline
\end{tabular}
\caption{\label{tab:nested_times}Times for the approximation of the Dirichlet-to-Neumann operator
 for the \textit{Laplace} and \textit{Helmholtz IV} problems via the accelerated nested dissection method.  }
\end{table}

\subsection{Performance with body loads}
\label{sec:bodyload}

In this experiment, the \textit{Random Laplacian I} problem is solved in a situation
with a non-zero body load $f$.  We assume however that the body load is restricted
to a small number $N_{\rm body}$ of nodes in the interior.
The locations of these nodes is assumed to be fixed.
Our objective is now to construct a solution
operator that constructs the vector of fluxes $\vct{v}$ on the boundary (the
discrete ``Neumann data'') given a vector $\vct{g}$ of Dirichlet data, and a
vector $\widehat{\vct{f}} \in \mathbb{R}^{N_{\rm body}}$ of body loads at the pre-scribed nodes.
This solution operator has two terms as follows:
\begin{equation}
\label{eq:general_solution_op}
\begin{array}{cccccccccccccc}
\vct{v} &=& \mtx{G}&\vct{g} &+& \mtx{F}&\widehat{\vct{f}},\\
N_{\rm boundary} \times 1 && N_{\rm boundary} \times N_{\rm boundary} & N_{\rm boundary} \times 1 && N_{\rm boundary} \times N_{\rm body} & N_{\rm body} \times 1
\end{array}
\end{equation}
where $N_{\rm boundary} = 4(n-1)$ denotes the number of points on the boundary.
The matrix $\mtx{G}$ is our by now familiar discrete Dirichlet-to-Neumann operator;
it is constructed in HBS form.
The matrix $\mtx{F}$ is a new solution operator that maps the interior body
load to boundary fluxes. Since $N_{\rm body}$ is small, the matrix $\mtx{F}$
is built in uncompressed form, and then approximated by a low-rank factorization.

Table \ref{tab:bodyloads} reports the computational times required to build
the solution operators $\mtx{F}$ and $\mtx{G}$ in (\ref{eq:general_solution_op})
for several different grid sizes and different values of $N_{\rm body}$. The table
also reports the relative errors when $N_{\rm body}$ random body loads are placed in a localized area inside the
domain.  Notice that for a small number of body loads the cost is close
to that of solving a pure boundary value problem.  As expected the computational
cost grows as the number of body loads is increased.


\begin{table}[ht]
 \centering
\begin{tabular}{|c|c|c|c|c|c|}
\hline
$N$ & $N_{\rm body}$ &  $T_{\rm build}$ 	&  $T_{\rm solve}$ &  $\epsilon_{\rm rel}$\\ \hline
\multirow{3}{*}{$512^2$} & 10 & 13.09 & 0.013 & 1.02e-6 \\
			  & 100 & 13.25 & 0.013 & 4.55e-7 \\
			& 1000 & 43.21 & 0.015 & 3.42e-7 \\ \hline
\multirow{3}{*}{$1024^2$} & 10 & 47.33 & 0.027 & 1.23e-6 \\
		      & 100 & 48.89 & 0.027& 6.46e-7 \\
		      & 1000 & 163.05 & 0.029 & 4.35e-7\\ \hline
\multirow{3}{*}{$2048^2$} & 10 & 256.57 & 0.55 & - \\ 
& 100 & 268.27& 0.58 & - \\
 & 1000 & 713.55 & 0.059 & - \\ \hline
\end{tabular}


\caption{\label{tab:bodyloads} Times for building the solution operators and applying the Dirichlet-to-Neumann operator when $N_{\rm body}$ body
loads are radomnly distributed in the domain.  The relative error $\epsilon_{\rm rel}$in the solution is also reported. }
\end{table}

\begin{remark}
In this section, we considered a special case where the body load $\vct{f}$ is restricted
to a small number of internal nodes. For the general case where the body load $\vct{f}$ is
supported on the entire domain, solution operators like (\ref{eq:general_solution_op}) can
still be constructed. In this case, the matrix $\mtx{F}$ is of size $N_{\rm boundary} \times N$,
and should be constructed in a data-sparse format analogous to the HBS format. This operator
can be both built and applied in $O(N)$ operations.
\end{remark}

\section{Conclusions and generalizations}
\label{sec:conc}
This paper presents a fast method for constructing the Dirichlet-to-Neumann operator
for elliptic problems with no body loads.
Numerical results indicate that the method scales linearly with the number of
discretization points $N$ for a variety of problems.
Since application of the solution operator scales linearly with the number of boundary
points (typically $O(N^{1/2})$), constructing the solution for multiple right-hand
sides is essentially free once the Dirichlet-to-Neumann operator is built.
For a problem involving approximately $16$ million unknowns, it takes about
$8$ minutes to build the solution operator, and $0.1$ seconds to apply it to a right-hand side.

The fast direct solver described here relies on the intermediate dense matrices being
compressible in the sense of being either of low rank, or having the HBS structure described in
Section \ref{sec:HBS}. It is currently not well understood exactly when this holds,
but the numerical experiments in Section \ref{sec:num} indicate that the property is
remarkably stable across a broad range of test problems.

In the interest of concision, this paper considered only an operator discretized by
a five-point stencil on a regular square grid. However, the scheme does not inherently
depend on the special form of either the stencil or the grid. We expect that the
generalization to other domains and other discretizations in 2D should in principle
be unproblematic, as long as the computational stencil is not too large. (Fast
construction of LU-decompositions of the stiffness operator on somewhat general
grids is reported in \cite{2010_ying_nesteddissection}.)

The scheme can also be generalized to problems in three dimensions; the simplistic
implementation described here would have $O(N^{3/2})$ complexity for the build stage,
and $O(N^{2/3})$ complexity for the solve stage. Given that classical nested dissection
in 3D has complexity $O(N^{2})$ and $O(N^{4/3})$ for the build and solve stages, this is
a substantial gain, especially for the solve stage. In principle, one could build a scheme
that uses accelerated matrix algebra internally inside the HBS representation to attain $O(N)$
complexity, but this would require significant work beyond that described in this paper.

\lsp

\noindent
\section*{Acknowledgements:} The work reported was supported by
NSF grants DMS0748488 and DMS0941476.

\bibliographystyle{amsplain}
\bibliography{IMANUM-refs}

\end{document}